%
\magnification=\magstep1
\input amstex
\UseAMSsymbols
\input pictex
\vsize=23truecm
\NoBlackBoxes
\parindent=18pt

   \font\rmk=cmr8    \font\itk=cmti8  \font\ttk=cmtt8

   \def\op{\text{\rm op}}
   \def\mod{\operatorname{mod}}

   \def\Hom{\operatorname{Hom}}
   \def\ind{\operatorname{ind}}
   \def\End{\operatorname{End}}
   \def\Ext{\operatorname{Ext}}

   \def\add{\operatorname{add}}
   
   \def\Cok{\operatorname{Cok}}

     \def\arr#1#2{\arrow <1.5mm> [0.25,0.75] from #1 to #2}

	\vglue1truecm
	\plainfootnote{}
	{\rmk 2010 \itk Mathematics Subject Classification. \rmk
	Primary 16G20, 16E65, 
			Secondary:
	16D90, 
       16E05, 
	16G10, 
	16G50. 
       16G70, 
       18G25. 
											 }

											 \centerline{\bf From submodule categories to preprojective algebras.}
											                      \bigskip
													      \centerline{Claus Michael Ringel and Pu Zhang}
													      				           \bigskip\medskip
																		   
{\narrower\narrower \rmk
Abstract: Let $\ssize S(n)$ be the category of invariant subspaces of nilpotent operators with nilpotency index at most $\ssize n$. Such submodule categories have been studied already in 1934 by Birkhoff, they have attracted a lot of attention in recent years, for example in connection with some weighted projective lines (Kussin, Lenzing, Meltzer). On the other hand, we consider the preprojective algebra $\ssize \Pi_n$ of type $\ssize \Bbb A_n$; the preprojective algebras
were introduced by Gelfand and Ponomarev, they are now of great interest, for example they form an important tool to study quantum groups (Lusztig) or cluster algebras (Geiss, Leclerc, 
Schr\"oer).

We are going to discuss the 
connection between the submodule category $\ssize \Cal S(n)$ and the module category $\ssize \mod \Pi_{n-1}$ 
of the preprojective algebra
$\ssize \Pi_{n-1}$. Dense functors $\ssize \Cal S(n) \to \mod \Pi_{n-1}$ are known to exist: one has 
 been constructed quite a long time ago by Auslander and Reiten, recently another one
by Li and Zhang. We will show that these two functors are full, objective functors
with index $\ssize 2n$, thus $\ssize \mod \Pi_{n-1}$ 
is obtained from $\ssize \Cal S(n)$ by factoring out an ideal which is generated by 
$\ssize 2n$ indecomposable objects. 

As a byproduct we also obtain new examples of
ideals in triangulated categories, namely ideals  $\ssize \Cal I$ in a triangulated
 category  $\ssize \Cal T$ which are 
generated by an idempotent such that the factor category $\ssize \Cal T/\Cal I$ 
is an abelian category.
\par}

      \bigskip
{\bf 1. Introduction.}
      \medskip
Let $k$ be a field. Let $S(n)$ be the category of invariant subspaces of nilpotent operators with nilpotency index at most $n$. For a detailed analysis of this category we refer to [RS2].
Let $k[x]$ be the polynomial ring in one variable $x$ with coefficients in $k$ and
$\Lambda_n = k[x]/\langle x^n\rangle$. The objects
of $\Cal S(n)$ are the pairs $(X,Y)$ where $Y$ is a $\Lambda_n$-module
and $X$ is a submodule of $Y$ (or the corresponding inclusion maps $u\:X \to Y$). 
We denote the indecomposable $\Lambda_n$-module
of length $i$ by $[i]$, and $[0]$ will denote the zero module.

Let  $\Pi_n$
be the preprojective algebra of type $\Bbb A_n$ and $\mod \Pi_n$ the category of the $\Pi_n$-modules of finite length.
The aim of this note is to show that $\mod \Pi_{n-1}$ 
is obtained from $\Cal S(n)$ by factoring out an ideal $\Cal I$ which is generated by 
$2n$ indecomposable objects --- actually, we will exhibit two possible choices for $\Cal I$.	
     \medskip 
Given an additive category $\Cal A$ and an ideal $\Cal I$ in $\Cal A$,
we denote by $\Cal A/\Cal I$ the corresponding factor category (it has the same
objects, and the homomorphisms in $\Cal A/\Cal I$ are the residue classes of
the homomorphisms in $\Cal A$ modulo $\Cal I$).
If $\Cal K$ is a class of objects of the category $\Cal A$,
we denote by $\langle \Cal K\rangle$ the ideal of $\Cal A$ given by all maps
which factor through a direct sum of objects in $\Cal K.$ Instead of writing
$\Cal A/\langle \Cal K\rangle$, we just will write $\Cal A/\Cal K.$
      \medskip
{\bf Theorem 1.} {\it Let $\Cal U$ be the set of objects of $\Cal S(n)$ 
which are of the form $([i],[j])$ with $i=j$ or $j = n$.
Let $\Cal V$ be the set of objects of
$\Cal S(n)$ which are of the form $([i],[j])$ with $i=j$ or $i = 0$.
Then the categories $\Cal S(n)/\Cal U$ and $\Cal S(n)/\Cal V$
are equivalent to $\mod\Pi_{n-1} $.

In particular, these categories  
are abelian categories with enough projective objects.
The indecomposable projective objects in $\Cal S(n)/\Cal U$
are the objects 
of the form $([0],[j])$ with $1\le j\le n-1$, 
the indecomposable projective objects in $\Cal S(n)/\Cal V$
are the objects of the form $([i],[n])$ with $1\le i\le n-1$.}
    \bigskip
If $\Cal A, \Cal B$ are Krull-Remak-Schmidt
categories, then we say that a functor $F\:\Cal A \to \Cal B$ is {\it objective,} provided its kernel 
is generated by identity maps of objects (see the appendix).
If $F$ is a dense, objective functor, then the number of isomorphism classes
of indecomposable objects $A$ in $\Cal A$ with $F(A) = 0$ is called its index. 
Thus Theorem 1 asserts that there are full, dense, objective functors 
$\Cal S(n) \to \mod\Pi_{n-1}$ with index $2n.$ These functors are the main
target of our considerations. Actually, the two functors
$F,G\:\Cal S(n) \to \mod\Pi_{n-1}$ which we will deal with have been exhibited
before: one of them has been constructed quite a long time ago by Auslander 
and Reiten [AR2], the other one recently by Li and Zhang [LZ], both are based
on general considerations by Auslander [A], published already in 1965.
What remains to be done is to show that the functors are full and objective.
     \medskip 
Both functors $F,G$ are given as compositions of functors which involve some
further module categories and also these intermediate categories seem to be
of interest. 

First of all, let $T_2(\Lambda_n)$ be the ring of upper triangular
$(2\times 2)$-matrices with coefficients in $\Lambda_n$. It is well-known (and easy
to see) that the category $\mod T_2(\Lambda_n)$-modules can be identified with the
category of maps between $\Lambda_n$-modules. Since the category $\Cal S(n)$ can
be considered as the category of monomorphisms of $\Lambda_n$-modules, we see
that $\Cal S(n)$ is a full subcategory of $\mod T_2(\Lambda_n)$, say with inclusion
functor $\iota$. There is a second full embedding $\epsilon\: \Cal S(n) \to 
\mod T_2(\Lambda_n)$, it sends the pair $(X,Y)$ in $\Cal S(n)$ to the canonical
projection $Y \to Y/X$.

Next, given an algebra $\Lambda$ of finite representation type, we denote by
$A(\Lambda)$ its (basic) Auslander algebra; it is defined as follows: let $E$ be a minimal
Auslander generator for $\Lambda$, this is the direct sum 
of all indecomposable $\Lambda$-modules, one from each isomorphism class, and
$A(\Lambda) = \End(E)^{\op}.$ 
Since the algebras $\Lambda_n$ are of
finite representation type, we may consider $A_n = A(\Lambda_n)$. Of course, in this
case $E = \bigoplus_{i=1}^n [i].$ 
In section 2, we study the
Auslander algebra $A_n$ and the category $\Cal F(n)$ of the torsionless
$A_n$-modules. Any Auslander algebra is quasi-hereditary, the Auslander algebras
$A_n$ are quasi-hereditary in a unique way and the category $\Cal F(n)$ is just
the category of $A_n$-modules with a $\Delta$-filtration.

The essential tool is the functor 
$$
 \alpha = \Cok\Hom_\Lambda(E,-) \: \mod T_2(\Lambda) \longrightarrow \mod A(\Lambda),
$$
where $\Lambda$ is of finite representation type and $E$ is a minimal Auslander 
generator $\Lambda$; note that $\alpha$ sends a morphism $f$ in $\mod \Lambda$ 
(thus the object $f$ of the category
$\mod T_2(\Lambda)$) to the cokernel of the induced map $\Hom_\Lambda(E,f)$; of course, 
$\Hom_\Lambda(E,f)$ is a map of $A(\Lambda)$-modules. This functor was considered already in 1965 by Auslander [A], and later by Auslander and Reiten in [AR1] and [AR2].  
Section 3 and parts of section 6 are devoted to this functor.

Finally, let us note that $\Pi_{n-1}$ is a factor ring of $A_n$, namely 
$\Pi_{n-1} = A_n/\langle e\rangle$, where $e$ is an idempotent of $A_n$ such that
$A_ne$ is indecomposable projective-injective. We will consider the functor
$$
 \delta\:\mod A_n \to \mod \Pi_{n-1}
$$
which sends any $A_n$-module $M$ to its largest factor module which is a $\Pi_{n-1}$-module,
thus to $M/A_neM.$ Properties of this functor $\delta$ will be discussed in section 5.
In particular, following [DR], we will look at the restriction $F_2$ of $\delta$ to the
subcategory $\Cal F(n).$

Altogether, we will deal with the following functors
$$
\hbox{\beginpicture
\setcoordinatesystem units <1.5cm,.7cm>
\put{$\Cal S(n)$} at 0 0
\put{$\mod T_2(\Lambda_n)$} at 2 0
\put{$\mod A_n$} at 4 0
\put{$\mod \Pi_{n-1}$} at 6 0
\arr{0.4 0.2}{1.3 0.2}
\arr{0.4 -.2}{1.3 -.2}
\arr{2.7 0}{3.5 0}
\arr{4.5 0}{5.4 0}
\put{$\iota$\strut} at .8 0.5
\put{$\epsilon$\strut} at .8 -.5
\put{$\alpha$\strut} at 3.1 0.3
\put{$\delta$\strut} at 4.9 0.3
\endpicture}
$$
The upper composition $F = \delta\alpha\iota$ is the functor studied by Li and Zhang [LZ],
the lower composition $G = \delta\alpha\epsilon$ is the one considered by
Auslander and Reiten [AR2]. We will see that {\it $F$ is a dense functor with kernel the
ideal $\langle \Cal U\rangle,$ and $G$ is a dense functor with kernel the
ideal $\langle \Cal V\rangle.$} This yields the first part of Theorem 1.

The image of the functor $\alpha\iota$ is precisely the
subcategory $\Cal F(n)$, thus we can write $F = F_2F_1$, where $F_1\:\Cal S(n) \to \Cal F(n)$
is the functor with $F_1(u) = \alpha(u)$, for $u\in \Cal S(n)$. 
It is known from the literature that the functors $F_1,F_2,G$ all are dense; our
main concern is to show that they are full and objective, with index $n,n,2n$, respectively. 
This will be done in sections 4, 5, 6, respectively. 

     \bigskip
In order to compare the functors $F$ and $G$, we have to take into account the
stable module category $\underline{\mod}\ \Pi_{n-1}$, it is obtained from the
module category $\mod \Pi_{n-1}$ by factoring out the ideal generated by the 
identity maps of the projective modules.
Since the algebra 
$\Pi_{n-1}$ is self-injective, the stable module category $\underline{\mod}\ \Pi_{n-1}$
is a triangulated category. We denote by
$$
 \pi\:\mod \Pi_{n-1} \longrightarrow  \underline{\mod}\ \Pi_{n-1}
$$
the canonical projection. Note that this is a full, dense, objective functor of index $n-1$,
its kernel is generated by the indecomposable projective $\Pi_{n-1}$-modules. 
    \bigskip 
We denote by $\Omega\:\underline{\mod}\ \Pi_{n-1} \to \underline{\mod}\ \Pi_{n-1}$
the syzygy functor, thus $\Omega(M)$ is the kernel of a projective cover of 
the $\Pi_{n-1}$-module $M$.
    \medskip
    {\bf Theorem 2.} {\it We have }
    $$
     \pi F = \Omega\pi G.
     $$
     \bigskip 
Section 8 draws the attention to the fact that in this setting we also obtain full,
dense, objective functors $\Cal T \to \Cal A$ with finite index,
such that $\Cal T$ is a triangulated category, $\Cal A$ an abelian category,
thus with ideals in triangulated categories which are generated by an idempotent
such that the corresponding factor categories are abelian.

In section 9 we provide illustrations
concerning the change of the Auslander-Reiten quiver of $\Cal S(n)$ when we factor out
the various ideals mentioned above.
	\bigskip
{\bf Acknowledgment.} The authors are indebted to the referee for a careful reading of
the paper and for many valuable comments concerning possible improvements of the paper, in particular, for suggesting to add section 9.

    \bigskip\bigskip
    {\bf 2. The Auslander algebra $A_n = A(\Lambda_n)$ and the subcategory $\Cal F(n)$ of
$\mod A_n$.}
     \medskip 
We recall that for $\Lambda$ a ring of finite representation type, $A(\Lambda)$ denotes
the basic Auslander algebra of $\Lambda$, it is the opposite of the
endomorphism ring of a minimal Auslander generator for $\Lambda$.
We consider here the special case of $\Lambda = \Lambda_n =
k[x]/\langle x^n\rangle$ and its Auslander algebra $A_n = (\End E)^{\op}.$
Note that $E = \bigoplus_{i=1}^n[i]$, where $[i]$ is the indecomposable
$\Lambda_n$-module of length $i$.

Let $P(i) = \Hom_\Lambda(E,[i]).$ This is an indecomposable projective $A_n$-module.
The inclusions $[i] \to [i+1]$ in the category $\mod \Lambda_n$ yield
a chain of inclusions
$$
 P(1) \subset P(2) \subset \cdots \subset P(n-1) \subset P(n).
 $$
 Let $\Delta(i) = P(i)/P(i-1)$ (with $P(0) = 0$). Note that $A_n$ is
 quasi-hereditary (for this ordering and only for this ordering) and the
 modules $\Delta(i)$ are the standard modules (but observe that the
 labeling of the simple $A_n$-modules exhibited here is the opposite of the
 labeling commonly used (see for example [DR]): in the present paper,
 it is the module $\Delta(1)$ which is projective
 and not $\Delta(n)$, and correspondingly, it is $P(n) = I(n)$ which is
 projective-injective and not $P(1)$).
 		      \medskip
Let $T(i) = P(n)/P(i-1)$ for $1\le i \le n$. Note that $T(i)$ is also the largest
submodule of $P(n-i+1)$ which is generated by $T(n)$. 
Let $T = \bigoplus_i T(i),$
this is the characteristic tilting module for the quasi-hereditary algebra
$A_n$.
		             \bigskip
			     {\bf Proposition 1.} {\it The following conditions are equivalent for an
			     $A_n$-module $M$.
			     \item{\rm (i)} $M$ is torsionless.
			     \item{\rm (ii)} $\Ext^1(M,T) = 0.$
			     \item{\rm (iii)} $M$ has a $\Delta$-filtration.
			     \item{\rm (iv)} The projective dimension of $M$ is at most $1$.
			     \item{\rm (v)} The injective envelope of $M$ is projective. \par}
			     	         \medskip
					 {\bf Remark 1.} We have avoided to refer to the labeling of the simple modules.
					 If we use our labels, so that $P(n) = I(n)$ is the unique indecomposable module which
					 is both projective and injective, then (v) can be reformulated as saying
					 that (vi) {\it the injective envelope $IM$ of $M$ 
is a direct sum of copies of $I(n)$,} or also that
(vii) {\it the socle of $M$ is a direct sum of copies of $S(n)$.}
      \medskip
      {\bf Remark 2.}
      It has been shown in [R] that any Auslander algebra is left strongly quasi-hereditary.
      This means that any module with a $\Delta$-filtration has projective dimension at
      most $1$ (the implication (iii) $\implies$ (iv)). 
Note that most of the relevant properties of  
strongly quasi-hereditary algebras have been considered already in [DR] and Proposition 1
is just a reminder.  
           \medskip
	   We denote by $\Cal F(n)$ the full subcategory of $\mod A_n$ given by the
	   $A_n$-modules which satisfy the equivalent conditions of Proposition 1. It follows
	   directly from (i) or also (v) that $\Cal F(n)$ is closed under submodules.
	   	     \medskip
		     Proof of Proposition 1.
		     (i) $\implies$ (ii). Let $M$ be torsionless. There is an embedding $M \to P$ with $P$
		     projective. Since the injective dimension of $T$ is $1$, the canonical map
		     $\Ext^1(P,T) \to \Ext^1(M,T)$ is surjective. But $\Ext^1(P,T) = 0$, thus
		     $\Ext^1(M,T) = 0.$

		     (ii) $\implies$ (iii). We have to show that
		     $\Ext^i(M,T) = 0$ for all $i\ge 1.$ Since $T$ has injective dimension $1$,
		     we only have to look at $i = 1$, but this is assertion (ii).
 
		     (iii) $\implies$ (iv). The modules $\Delta(i)$ have projective dimension at
		     most 1.

		     (iv) $\implies$ (v). Let $0 @>>> P_1 @>u>> P_0 @>>> M @>>> 0$ be
		     a  projective resolution. We have to show that $M$ embeds into a module
		     which is both injective and projective.
		     Consider the injective envelopes of
		     $v_i\:P_i \to I(P_i)$ for $i=0,1$,
		     this yields a commutative diagram with exact rows
		     $$
		     \CD
		      0 @>>> P_1  @>u>> P_0  @>>> M  @>>> 0 \cr
		       @.    @Vv_1 VV      @VVv_0 V      @VVf V \cr
		        0 @>>> I(P_1)  @>u'>> I(P_0)  @>>> M'  @>>> 0 .
			\endCD
			$$
			Since $v_0$ is injective, the snake lemma yields an embedding of
			the kernel $K$ of $f$ into the cokernel of $v_1$. For any projective
			module $P$ with injective envelope $IP$,
the cokernel $IP/P$ embeds into a projective-injective module
(since the dominant dimension of $\Lambda_n$ is at least 2).

On the other hand, $u'$ is a monomorphism and $I(P_1)$ is injective,
thus $u'$ is a split monomorphism. Thus $M'$ is a direct summand of
$I(P_0)$ and this module is projective-injective.
The exact sequence $0 \to K \to M \to M'$ shows that $M$ embeds into
$IK\oplus M'$ which is projective-injective.

(v) $\implies$ (i). If the injective envelope of $M$ is projective, then $M$
embeds into a projective module, thus $M$ is torsionless.
\hfill $\square$
       \bigskip
       {\bf Proposition 2.} {\it If $M$ is an $A_n$-module of projective dimension at most $1$
       and generated by $P(n)$, then $M$ belongs to $\add T$.}
           \medskip
	   Proof: Write $M = P/U$ where $P$ is a direct sum of copies of $P(n)$ and
	   $U$ is a submodule of $P$. Since the projective dimension of $M$ is at most 1,
	   the submodule $U$ is projective. Now $P$ is injective, thus we may embed
	   an injective envelope $IU$ of $U$ into $P$, this is a direct summand of $P$,
	   say $P = IU\oplus C$ for some submodule $C$. Since $P$ is a direct sum of copies
	   of $P(n)$, also $C$ is a direct sum of copies of $P(n) = T(1).$ On the other
	   hand, the exact sequences $0 \to P(i) \to P(n) \to T(i+1) \to 0$ for $1\le i \le n$
	   (with $T(n+1) = 0$) show that for any indecomposable projective module $P'$,
	   the module $I(P')/P'$ belongs to $\add T$, thus $IU/U$ belongs to $\add T.$
	   Altogether we see that $P/U = IU/U\oplus C$ belongs to $\add T.$
\hfill$\square$
	   	         \medskip
			 {\bf Remark.} The modules generated by $P(n) = T(1)$ are the modules
			 with a $\nabla$-filtration (see [DR]), thus proposition 2 just asserts that modules
			 which have both a $\Delta$-filtration and a $\nabla$-filtration belong to
			 $\add T.$
			       \bigskip\bigskip
			       
{\bf 3. The functor $\alpha\:\mod T_2(\Lambda) \to \mod A(\Lambda)$.}
			            \medskip	       
Here we deal with an arbitrary representation-finite algebra $\Lambda$ 
with minimal Auslander generator $E$ and consider  the corresponding Auslander algebra
$A(\Lambda) = (\End E)^{\op}$.
We have 
denoted $T_2(\Lambda)$  the ring of all upper triangular $(2\times 2)$-matrices with coefficients in $\Lambda$. The category of $T_2(\Lambda)$-modules may be seen as
the category of all morphisms $X \to Y$ in $\mod\Lambda$. We consider the functor
$$
  \alpha\:\mod T_2(\Lambda) \to \mod A(\Lambda)
$$
defined by $\alpha(f) = \Cok\Hom_\Lambda(E,f)$ for a morphism $f$ in $\mod\Lambda$
    \medskip
    {\bf Proposition 3.} {\it Let $\Cal X$ consist of the two objects
    $(1\:E \to E)$ and $(E \to 0)$ in $\mod T_2(\Lambda)$.
    The functor
    $\alpha\:\mod T_2(\Lambda) \to \mod A(\Lambda)$
        yields an equivalence
    $(\mod T_2(\Lambda))/\Cal X \to \mod A(\Lambda).$} 
    
\medskip 
Thus, if the number of isomorphism classes of indecomposable $\Lambda$-modules
is $m$, {\it then $\alpha\:\mod T_2(\Lambda) \to \mod A(\Lambda)$ is a full, dense, objective functor
with index $2m$.}
       \medskip
       For a slightly weaker statement we refer to Theorem 1.1 in [AR2].
       Let us recall and complete the proof. Given a morphism $f\:X \to Y$ in $\mod\Lambda$,
       we obtain an exact sequence
       $$
        \Hom_\Lambda(E,X) @>\Hom_\Lambda(E,f)>> \Hom_\Lambda(E,Y) @>>> \alpha(f) @>>> 0.
	$$
	Since the $A(\Lambda)$-modules $\Hom_\Lambda(E,X)$ 
and $\Hom_\Lambda(E,Y)$ are projective, we
	obtain in this way a projective presentation of $\alpha(f).$ Conversely, given
	an $A(\Lambda)$-module $M$, take a projective resolution
	$P_1 @>p>> P_0 \to M \to 0.$ Now, the category of projective $A(\Lambda)$-modules is
	equivalent to the category $\mod \Lambda$, thus we can assume that
	there are $\Lambda$-modules $X$ and $Y$ such that
	$P_1 = \Hom_\Lambda(E,X)$ and $P_0 = \Hom_\Lambda(E,Y)$ and a $\Lambda$-homomorphism $f\:X \to Y$
	such that $\Hom_\Lambda(E,f) = p.$ In this way, we see that the functor is dense.
	Similarly, starting with an $A(\Lambda)$-homomorphism $M \to M'$, we can lift it
	to projective presentations of $M$ and $M'$ and using again the fact that
	the category of projective $A(\Lambda)$-modules is
	equivalent to the category $\mod \Lambda$, we see that $M \to M'$ is in the image
	of the functor $\alpha.$ Thus, it remains to calculate the kernel of $\alpha.$

	Of course, under this
	functor $\alpha$ the two objects in $\Cal X$ are sent to zero. Thus the ideal
	$\langle \Cal X\rangle$ is contained in the kernel of $\alpha$.
	Conversely, assume that there is given a map
	$$
	 (g_1,g_0)\:(f\:X_1\to X_0) \longrightarrow (f'\:X_1' \to X_0')
	 $$
	 (thus $g_0f = f'g_1$),
	 such that $\alpha(g_1,g_0) = 0.$ Thus the
	 following diagram commutes and its rows are projective presentations:
	 $$
	 \CD
	  \Hom_\Lambda(E,X_1) @>\Hom_\Lambda(E,f)>> \Hom_\Lambda(E,X_0) @>e>> \alpha(f) @>>> 0\cr
	    @V\Hom_\Lambda(E,g_1)VV    @VV\Hom_\Lambda(E,g_0)V    @VV{\alpha(g_1,g_0)=0}V \cr
	     \Hom_\Lambda(E,X_1') @>\Hom_\Lambda(E,f')>> \Hom_\Lambda(E,X_0') @>e'>> \alpha(f') @>>> 0.
	     \endCD
	     $$
	     Let us show that the map $(g_1,g_0)$ factors through $[1\ 0]\:X_0\oplus X_1 \to X_0.$
	     Now $e'\Hom_\Lambda(E,g_0) = 0$, thus there is a map $\widetilde h\:\Hom_\Lambda(E,X_0) \to
	     \Hom_\Lambda(E,X_1')$ such that $\Hom_\Lambda(E,f')\widetilde h = \Hom_\Lambda(E,g_0).$ Again using the
	     equivalence of $\mod \Lambda_n$ and the category of projective $A_n$-modules,
	     we get a map $h\:X_0 \to X_1'$ with $\widetilde h = \Hom_\Lambda(E,h)$ and $g_0 = f'h.$
	     Therefore $f'hf = g_0f = f'g_1$, and thus $f'(g_1-hf) = 0.$ As a consequence, the following
	     diagram commutes
	     $$
	     \CD
	      X_1 @>f>> X_0  \cr
	        @V{\left[\smallmatrix f\cr 1\endsmallmatrix\right]}VV    @VV1V \cr
		X_0\oplus X_1 @>{\left[\smallmatrix 1&0\endsmallmatrix\right]}>> X_0  \cr
		  @V{\left[\smallmatrix h& g_1-hf\endsmallmatrix\right]}VV    @VVg_0V \cr
		   X'_1 @>f'>> X'_0
		   \endCD
		   $$
		   and the composition of the vertical maps is just $(g_1,g_0)$. Since $E$
		   is an Auslander generator, all $\Lambda_n$-modules belong to $\add E$, thus
		   $[1\ 0]\:X_0\oplus X_1 \to X_0$ belongs to $\add \Cal X.$
		   This shows that the map $(g_1,g_0)$ belongs to $\langle\Cal X\rangle.$
\hfill$\square$
		        \bigskip\bigskip
			{\bf 4. The functor $F_1:\Cal S(n) \to \Cal F(n)$.}
			     \medskip
			     It may be appropriate to
			     focus first the attention to the category $\Cal S(n)$ itself.
			     An object $(X,Y)$ of $\Cal S(n)$ with inclusion map $u\:X \to Y$
			     will also be denoted by $u$. This stresses the fact that objects of
			     $\Cal S(n)$ are given by
			     maps in the category $\mod \Lambda_n$, thus we consider $\Cal S(n)$ as a
			     full subcategory of the category of $T_2(\Lambda_n)$-modules. As we have mentioned,
we denote
by $\iota\:\Cal S(n) \to \mod T_2(\Lambda_n)$ the inclusion functor.
Note that $\Cal S(n)$ turns out to be just the category
of all Gorenstein-projective $T_2(\Lambda_n)$-modules, see for example [Z].
   \medskip
   We consider the restriction $\alpha\iota$ of
   $\alpha$ to the full subcategory $\Cal S(n)$ of $\mod T_2(\Lambda_n)$, thus
   $$
    \alpha\iota(u) = \Cok \Hom_\Lambda(E,u).
    $$
    for $(u\:X\to Y)$ in $\Cal S(n)$. Since $u$ is a monomorphism,
    also $\Hom_\Lambda(E,u)$ is a monomorphism, thus there is the following exact
    sequence
    $$
     0 \to \Hom_\Lambda(E,X) \to \Hom_\Lambda(E,Y) \to \alpha(u) \to 0.
     $$
     Since both $\Hom_\Lambda(E,X)$ and $\Hom_\Lambda(E,Y)$ are projective $A_n$-modules,
     we see that the projective dimension of $\alpha(u)$ is at most $1$,
     thus $\alpha(u)$ belongs to $\Cal F(n).$ 
This shows that we can consider $\alpha\iota$ as a functor $\Cal S(n) \to \Cal F(n)$,
we denote it by $F_1$ (thus $F_1(u) = \alpha(u)$, for $u\in \Cal S(n)$). 

     \medskip
     Under the functor $F_1$, we have
     $$
     \alignat 2
      F_1([i],[i]) &= 0     &&1\le i \le n, \cr
       F_1([i],[n]) &= P(n)/P(i) = T(i\!+\!1) &\qquad\text{for}\qquad &0\le i \le n\!-\!1, \cr
        F_1([0],[j]) &= P(j)  &&1\le j\le n\!-\!1. \cr
	\endalignat
	$$
		\medskip
		{\bf Proposition 4.} {\it Let $\Cal U_1$ be the set of objects of
		$\Cal S(n)$ of the form $([i],[i])$ with $1\le i \le n$.
		The functor $F_1$ yields an equivalence between
		the factor category of $\Cal S(n)/\Cal U_1$ and the category $\Cal F(n)$.
		It maps the pairs $([i],[n])$ with $0 \le i \le n\!-\!1$ to $T(i+1)$
		and the pairs $([0],[j])$ with $1\le j \le n\!-\!1$ to $P(j)$.}
		    \medskip
		    Proof. By definition, $F_1$ is the restriction of the functor $\alpha$
		    to the subcategory $\Cal S(n)$ of $\mod T_2(\Lambda_n)$. We have noted already
		    that the image of $F_1$ consists of modules of projective dimension at most $1$.
		    Also, conversely, if $M$ is a $A_n$-module of projective dimension at most $1$,
		    say with a projective presentation
		    $$
		     0 @>>> P_1 @>p>> P_0 @>>> M @>>> 0,
		     $$
		     then we can write $p = \Hom_\Lambda(E,f)$ for some map $f\:X_1 \to X_0$ in $\mod \Lambda_n$.
		     Clearly, $f$ has to be a monomorphism, since $E$ is an Auslander generator.
		     Thus we can assume that $f$ belongs to $S(n)$ and we have $F_1(f) = M.$

		     According to Proposition 3 the kernel of the functor $F_1$ is given by all
		     morphisms in $\Cal S(n)$ which factor through a $T_2(\Lambda_n)$-module of the
		     form $[1\ 0]\:V \oplus V' \to V.$ Thus assume we have the following commutative
		     diagram
		     $$
		     \CD
		      X @>u>> Y  \cr
		        @V{\left[\smallmatrix g\cr g'\endsmallmatrix\right]}VV    @VVg''V \cr
			 V\oplus V' @>{\left[\smallmatrix 1&0\endsmallmatrix\right]}>> V  \cr
			   @V{\left[\smallmatrix h& h'\endsmallmatrix\right]}VV    @VVh''V \cr
			    X' @>u'>> Y'.
			    \endCD \tag$*$
			    $$
			    The commutativity of the lower square means that we have
			    $$
			     [h'',0] = h''[1,0] = u'[h,h'] = [u'h,u'h'],
			     $$
			     thus $u'h' = 0.$ Since $u'$ is a monomorphism, it follows that $h' = 0.$ But then
			     also the diagram
			     $$
			     \CD
			      X @>u>> Y  \cr
			        @Vg VV    @VVg''V \cr
				 V @>1>> V  \cr
				   @Vh VV    @VVh''V \cr
				    X' @>u'>> Y'.
				    \endCD
				    $$
				    commutes and the composition of the vertical maps is the same as the composition
				    of the vertical maps in $(*)$. This shows that the kernel of the functor $F_1$
				    consists of all the morphisms in $S(n)$ which factor through objects of the form
				    $(1\:V \to V)$, but this is just the ideal $\Cal U_1.$
\hfill$\square$
				           \bigskip\bigskip
   {\bf  5. The functor $\delta\:\mod A_n \to \mod \Pi_{n-1}$ 
and its restriction $F_2$ to $\Cal F(n)$.}
      \medskip
      Recall that the indecomposable
      projective $A_n$-modules are $P(i) = \Hom_\Lambda(E,[i])$ with $1\le i \le n$,
      if necessary, we will denote them by $P_A(i)$. The module $P(n)$ is also injective
and we may choose an idempotent $e(n)$ in $A_n$ such that $P(n) = A_ne(n)$. 
Let $\Pi_n$ be the preprojective algebra of type $\Bbb A_n$; note that
$A_n/\langle e(n)\rangle = \Pi_{n-1}$ (see [DR], Theorem 3, Theorem 4 and Chapter 7). 
	          \medskip
		  We consider in $\mod A_n$ the following torsion pair: the torsion modules
		  are the modules generated by $P(n)$, thus the modules with top being a direct
		  sum of copies of $S(n)$, the torsionfree modules are the modules which do not
		  have $S(n)$ as a composition factor (thus the torsionfree modules form a
		  Serre subcategory). Note that the torsionfree modules are just the $A_n$-modules
		  $M$ with $e(n)M = 0$, thus these are the $\Pi_{n-1}$-modules.

		  Given an $A_n$-module $M$, let $tM$ be its torsion submodule, and
		  $$
		   \delta M = M/tM.
		   $$
		   The indecomposable
		   projective $\Pi_{n-1}$-modules are factor modules of the modules $P_A(i)$
		   with $1\le i \le n\!-\!1$, we denote them by $P_\Pi(i) = P_A(i)/tP_A(i).$

		   Since we want to keep track of the selected objects $([i],[j])$
		   in $\Cal W$, let us repeat what happens under the functor $\delta$:
		   $$
		   \alignat 2
		    \delta(T(i)) &= 0 &\qquad\text{for}\qquad &1\le i \le n, \cr
		     \delta(P_A(j)) &= P_\Pi(j)  &&1\le j\le n\!-\!1. \cr
		     \endalignat
		     $$
			\medskip
			{\bf Proposition 5.} {\it
			The functor $\delta$ yields an equivalence between
			the factor category $\Cal F(n)/T$ and the category $\mod \Pi_{n-1}$.
			It maps $P_A(j)$ with $1\le j \le n\!-\!1$ to $P_\Pi(j)$.}
			   \medskip
			   Proof. First, let us show that the functor $\delta$ is dense. Thus, let $N$
			   be a $\Pi_{n-1}$-module, but consider it as an $A_n$-module. Actually,
			   all the modules to be considered now are $A_n$-modules; in particular
			   $I(i)$ will denote (as before) the indecomposable injective $A_n$-module
			   corresponding to the vertex $i$.

			   Let $u\:N \to IN$
			   be an injective envelope of $N$ (as an $A_n$-module!)
			   and $p\:PIN \to IN$ a projective cover of $IN$.
			   Let $T'$ be the kernel of $p$. Now $IN$ is a direct sum of modules of the form
			   $I(i)$ with $1\le i \le n-1$ and the exact sequence $0 \to T(i+1) \to P(n) \to I(i) \to 0$
			   shows that the kernel of the projective cover $P(n) \to I(i)$ is $T(i+1)$.
			   Thus we see that $T'$ is a direct sum of copies of modules of the form $T(j)$
			   with $2\le j \le n.$

			   Forming the induced exact sequence with respect to $u$, we obtain the following
			   commutative diagram with exact rows:
			   $$
			   \CD
			    0 @>>> T' @>v'>> \widetilde N @>p'>> N @>>> 0 \cr
			     @.     @|       @VVu'V           @VVu V \cr
			      0 @>>> T' @>v>>  PIN     @>p>>   IN @>>> 0
			      \endCD
			      $$
			      The monomorphism $u'$ shows that $\widetilde N$ is torsionless, thus in
			      $\Cal F(n)$. Clearly, $tT' = T'$ and $tN =0$, thus $t\widetilde N$ is
			      the image of $v'$ and therefore $\delta(\widetilde N) = N.$ This shows that
			      $\delta$ is dense. In the same way, we see that $\delta$ is also full.
			      Namely, given a morphism $f\:N_1 \to N_2$ of $\Pi_{n-1}$-modules, we can extend
			      $f$ to a morphism $f'\:I(N_1) \to I(N_2)$ and then lift it to a morphism $f''\:PI(N_1)
			      \to  PI(N_2).$ Using the pullback property of the commutative square
			      $$
			      \CD
			        \widetilde N_2 @>p'>> N_2 \cr
				    @VVu'V           @VVu V \cr
				       PI(N_2)     @>p>>   I(N_2)
				       \endCD
				       $$
				       we finally obtain a map $\widetilde f\:\widetilde N_1 \to \widetilde N_2$
				       such that $\delta(\widetilde f) = f.$

				       It remains to determine the kernel of $\delta$. Since the modules $T(i)$ are
				       generated by $T(1)$, they are torsion modules, thus $\delta(T(i)) = 0$ for
				       all $1\le i \le n$. Conversely, let $M_1,M_2$ belong to $\Cal F(n)$ and let
				       $f\:M_1 \to M_2$ be a homomorphism
				       such that $\delta(f) = 0.$ This means that the image of $f$
				       is contained in the torsion submodule $tM_2.$ Now $tM_2$ is a submodule
				       of $M_2 \in \Cal F(n)$ and $\Cal F(n)$ is closed under submodules, thus
				       $tM_2$ has projective dimension at most $1$. On the other hand, $tM_2$ is generated
				       by $P(n)$. Thus Proposition 2 asserts that $tM_2$ is in $\add T$. It follows that
				       $f$ belongs to the ideal $\langle T\rangle.$

The last assertion has been mentioned already before.
\hfill$\square$
				           \medskip
					   {\bf Remark 1.} In order to show the density of the functor $\delta$,
					   the paper [DR] started with a $\Pi_{n-1}$-module $N$, considered it
					   as an $A_n$-module and used a universal extension
					   $$
					    0 \to T' \to \widetilde N \to N \to 0
					    $$
					    of $N$ by a module $T'$ in $\add T$. Of course, the pullback recipe
					    given in the proof above provides such a universal extension.
					          \medskip
						  {\bf Remark 2.} Looking at the exact sequence
						  $$
						   0 \to T' \to \widetilde N \to N \to 0
						   $$
						   constructed in the proof of Proposition 5,
						   one may decompose $T' = \bigoplus_i T(i)^{n(i)}$. Then $n(1) = 0$ and
						   {\it $n(i+1)$ is precisely the multiplicity of $S(i)$ in the socle of
						   $N$.}
							\bigskip
							Combining Propositions 4 and 5 we obtain the first part of the following proposition:
							\medskip
{\bf Proposition 6.} {\it The functor $F\:\Cal S(n) \to \mod\Pi_{n-1}$ is full, dense and
its kernel is $\langle\Cal U\rangle$, thus $F$ induces an equivalence 
$\Cal S(n)/\Cal U \to \mod\Pi_{n-1}.$

Under such an equivalence, the objects of $\Cal S(n)$
of the form $([0],[j])$ with $1\le j\le n-1$ correspond to the indecomposable
projective $\Pi_{n-1}$-modules.}
	   \medskip 
The last assertion relies on the
fact that the set of projective objects of an abelian category is uniquely determined.
Since  $F([0],[j]) = P_{\Pi}(j)$, for $1\le j \le n-1$,
we know that the objects $([0],[1]),\dots ,([0],[n-1])$ in $\Cal S(n)/\Cal U$ are
the indecomposable projective objects in $\Cal S(n)/\Cal U$. Thus any equivalence
between $\Cal S(n)/\Cal U$ and $\mod \Pi_{n-1}$ sends these objects
to the indecomposable projective $\Pi_{n-1}$-modules.
\hfill$\square$
     \bigskip\bigskip
     {\bf 6. The functor $G\:\Cal S(n) \to \mod\Pi_{n-1}$.}
          \medskip
In this section, we are going to 
analyze the functor $G\:\Cal S(n) \to \mod\Pi_{n-1}$.
The essential observation is due to Auslander and Reiten, it 
is valid in the general setting of dealing with a 
representation-finite algebra $\Lambda$ as discussed in section 3.
	\medskip
{\bf Lemma.} {\it Let $\Lambda$ be a representation-finite algebra with
Auslander generator $E$ and Auslander algebra $A = (\End E)^{\op}$.
Let $f$ be a morphism of $\Lambda$-modules. Then $f$ is an epimorphism if and
only if $\Hom_A(\Hom_\Lambda(E,\Lambda),\alpha(f)) = 0.$}
	\medskip
We may refer to Proposition 4.1 in [AR1], see also
the formulation stated in [AR2] just after the proof of Theorem 1.1.
But the reader should observe that for the proof one just has to use twice 
Yoneda isomorphisms. Namely, starting with the
map $f\:X \to Y$ of $\Lambda$-modules, consider the 
map $\Hom_\Lambda(E,f)\:\Hom_\Lambda(E,X) \to \Hom_\Lambda(E,Y)$  
of $A$-modules; by definition, $\alpha(f)$ is its cokernel,
thus there is the exact sequence
$$
\CD
  \Hom_\Lambda(E,X) @>\Hom_\Lambda(E,f)>> \Hom_\Lambda(E,Y) @>>> \alpha(f) @>>> 0.
\endCD
$$
If we apply the functor $H = \Hom_A(\Hom_\Lambda(E,\Lambda),-)$, we obtain the
sequence
$$
\CD
  H(\Hom_\Lambda(E,X)) @>H(\Hom_\Lambda(E,f))>> H(\Hom_\Lambda(E,Y)) @>>> H(\alpha(f)) @>>> 0,
\endCD
$$
it is exact, since $\Hom_\Lambda(E,\Lambda)$ is a projective $A$-module. In this way,
we see that 
$H(\alpha(f)) = \Hom_A(\Hom_\Lambda(E,\Lambda),\alpha(f))$ is 
the cokernel of $H(\Hom_\Lambda(E,f)).$ There is the following commutative diagram
$$
\CD 
  H(\Hom_\Lambda(E,X)) @>H(\Hom_\Lambda(E,f))>>  H(\Hom_\Lambda(E,Y)) \cr
      @VVV                                  @VVV \cr
  \Hom_\Lambda(\Lambda,X) @>\Hom_\Lambda(\Lambda,f)>> \Hom_\Lambda(\Lambda,Y)
\endCD
$$
where the vertical maps are Yoneda isomorphisms. The lower map may be identified with the map $f$
(again, this is due to Yoneda isomorphisms). 
Thus $f$ is an epimorphism if and only if $\Hom_\Lambda(\Lambda,f)$ is an
epimorphism if and only if $\Hom_A(\Hom_\Lambda(E,\Lambda),\Hom_\Lambda(E,f))$ 
is an epimorphism, if and only of the cokernel of  
$\Hom_A(\Hom_\Lambda(E,\Lambda),\Hom_\Lambda(E,f))$ is zero, thus if and only if 
$\Hom_A(\Hom_\Lambda(E,\Lambda),\alpha(f)) = 0$. \hfill $\square$

	\bigskip  
Let us return to the special case of $\Lambda = \Lambda_n$ and 
$A = A_n$. Recall that we have chosen an idempotent $e(n)$ in $A$
such that $A_ne(n) = P(n) = \Hom_{\Lambda_n}(E,[n]) = \Hom_{\Lambda_n}(E,\Lambda_n)$
and we have identified $A_n/\langle e(n)\rangle = \Pi_{n-1}$. Thus, given an
$A_n$-module $M$, the condition 
$\Hom_{A_n}(\Hom_ {\Lambda_n}(E,\Lambda_n),M) = 0$ can be rewritten as
$\Hom_{A_n}(P(n),M) = 0,$ thus as $e(n)M = 0,$ and the modules $M$ with
this property are just the $\Pi_{n-1}$ modules. Thus, in our case the Lemma
asserts the following: 
{\it Let $f$ be a morphism of $\Lambda_n$-modules. Then $f$ is an epimorphism if and
only if $\alpha(f)$ is a $\Pi_{n-1}$-module.}
	\medskip 
We consider now the functor $G = \delta\alpha\epsilon.$ 
If  $u\:X \to Y$ in $\Cal S(n)$, then $\epsilon(u)$ is an epimorphism,
thus $\alpha\epsilon(u)$ is a $\Pi_{n-1}$-module and therefore 
$$
   G(u) = \delta\alpha\epsilon(u) = \alpha\epsilon(u).
$$
This shows that $G\:\Cal S(n) \to \mod\Pi_{n-1}$ is defined by
$$
 G(X,Y) = \Cok\Hom_\Lambda(E,Y\to Y/X),
$$
we form the cokernel map $q = \epsilon(u)\:Y \to Y/X$ of the inclusion map
$u\:X\to Y$, apply the functor $\Hom_\Lambda(E,-)$ so that we obtain a map
$$
  \Hom_\Lambda(E,q)\:\Hom_\Lambda(E,Y) \to \Hom_\Lambda(E,Y/X)
$$
and take its cokernel $\Cok\Hom_\Lambda(E,q)$.

{\bf Remark.} We may phrase the definition of $G$ also differently:
  Let $\Cal Q(n)$ be the category of all epimorphisms $(Y \to Z)$ in $\mod \Lambda_n$,
  this is again a subcategory of $\mod T_2(\Lambda_n)$ and we may look at the restriction
  of the functor $\alpha$ to $\Cal Q(n)$. Of course, there is an obvious
  categorical equivalence $\Cal S(n) \to \Cal Q(n)$, it sends an object
  $(u\:X \to Y)$ to the canonical map $(Y \to Y/X),$ this is just $\epsilon(u)$ and 
$G = \alpha\epsilon.$

Instead of looking at the categories $\Cal S(n)$ and $\Cal Q(n)$,
we may consider the category $\Cal E(n)$ of all short exact sequences in the category
$\mod \Lambda_n$. There are forgetful functors $\Cal E(n) \to \Cal S(n)$
and $\Cal E(n) \to \Cal Q(n)$ which send an exact sequence
$0 @>>> X @>u>> Y @>q>> Z @>>> 0$ to $u$ or $q$, respectively. Obviously, both functors are
categorical equivalences.
	        \medskip
		Under the functor $G$, we have
		$$
		\alignat 2
		 G([i],[i]) &= 0     &&1\le i \le n, \cr
		  G([i],[n]) &= P_{\Pi}(n\!-\!i) &\qquad\text{for}\qquad &1\le i \le n\!-\!1, \cr
		   G([0],[j]) &= 0  &&1\le j\le n. \cr
		   \endalignat
		   $$
(The second assertion is seen as follows: Let $q\:[n] \to [n-i]$ be a cokernel map
for the inclusion map $u\:[i] \to [n]$. 
If we apply $\Hom_\Lambda(E,-)$ to the short exact sequence with maps $u$ and $q$, we obtain the
exact sequence 
$$
\CD
 0 \to P(i) @>\Hom_\Lambda(E,u)>> P(n) @>\Hom_\Lambda(E,q)>> P(n-i).
\endCD
$$   
Now the embedding $\Hom_\Lambda(E,u)$ of $P(i)$ into $P(n)$ has as cokernel the module $T(i+1)$.
In this way, $T(i+1)$ is embedded into $P(n-i)$; the image $U$ of this embedding is the
largest submodule of $P(n-i)$ generated by $P(n)$. This shows that $P(n-i)/U$ is
equal to $P_{\Pi}(n\!-\!i)$. On the other hand, $P(n-i)/U$ is just the cokernel of
$\Hom_\Lambda(E,q)$, thus equal to $G([i],[n]).$)
			\medskip
{\bf Proposition 7.} {\it The functor $G\:\Cal S(n) \to \mod\Pi_{n-1}$ is full, dense
and $\langle V\rangle$ is its kernel. Thus $G$ induces an equivalence
$\Cal S(n)/\Cal V \to \mod\Pi_{n-1}.$

Under such an equivalence, the objects of $\Cal S(n)$
of the form $([i],[n])$ with $1\le i\le n-1$ correspond to the indecomposable
projective $\Pi_{n-1}$-modules.}
	   \medskip
Proof.
Proposition 3 asserts that $\alpha$ is a full  functor from
$\Cal Q(n)$ onto the category $\mod \Pi_{n-1}$ and that
its kernel is the set of morphisms which factor through an object of the
form $([1,0]\:V\oplus V' \to V).$ Under the equivalence $\Cal S(n) \to
\Cal Q(n)$, the objects of the form $([1,0]\:V\oplus V' \to V)$ in $\Cal Q(n)$
correspond to the objects of the form $(\left[\smallmatrix 0\cr 1
\endsmallmatrix\right]\:V' \to V\oplus V')$ in $\Cal S(n)$,
but these are precisely the objects
in $\add \Cal V.$ This shows that $G= \alpha\epsilon$ yields an equivalence
$\Cal S(n)/\Cal V \to \mod\Pi_{n-1}$.

By proposition 3, we know that the functor $\alpha\:\mod T_2(\Lambda_n) \to \mod A_n$
is full and dense. Let us show that the restriction of $\alpha$ to $\Cal Q(n)$ is a dense
functor $\Cal Q(n) \to \mod \Pi_{n-1}$. If $M$ is a $\Pi_{n-1}$-module, the density
of $\alpha$ provides a map $f$ of $\Lambda_n$-modules such that $\alpha(f)$ is
isomorphic to $M$. But since $\alpha(f)$ is a $\Pi_{n-1}$-module, we know that
$f$ is an epimorphism (see the reformulation of the Lemma above), thus 
$f$ belongs to $\Cal Q(n)$.

The last assertion of Proposition 7 relies again on the
fact that the set of projective objects of an abelian category is uniquely determined
by the categorical structure. 
\hfill$\square$
	\medskip
This completes the proof of Theorem 1.
   \bigskip\bigskip

   {\bf 7. Comparison of the functors $F$ and $G$.}
     \medskip
     We want to compare the functor $F = F_2F_1$ and $G$. In order
     to prove the equality $\pi F = \Omega\pi G,$ we
     start with an object $u\:X \to Y$ in $\Cal S(n)$, thus with an exact sequence
     $$
      0 @>>> X @>u>> Y @>q>> Z @>>> 0.
      $$
      We apply the functor $\Hom_\Lambda(E,-)$ and obtain the exact sequence
      $$
       0 @>>> \Hom_\Lambda(E,X) @>\Hom_\Lambda(E,u)>> \Hom_\Lambda(E,Y) @>\Hom_\Lambda(E,q)>> \Hom_\Lambda(E,Z)
       $$
       Now, the cokernel of $\Hom_\Lambda(E,q)$ is $G(u)$.
       The cokernel of $\Hom_\Lambda(E,u)$ and thus the image of $\Hom_\Lambda(E,q)$ is
       $F_1(u)$, thus there is the following exact sequence
       $$
        0 @>>> F_1(u) @>>> \Hom_\Lambda(E,Z) @>>> G(u) @>>> 0
	$$
	and we can assume that the map $F_1(u) \to \Hom_\Lambda(E,Z)$ is an
	inclusion map.
	Since $G(u)$ is a $\Pi_{n-1}$-module, we have $tG(u) = 0$,
	therefore $t\Hom_\Lambda(E,Z) \subseteq F_1(u)$, and therefore $t\Hom_\Lambda(E,Z) =
	tF_1(u).$
	Thus we have the following
	exact sequence:
	$$
	 0 @>>> F_1(u)/tF_1(u) @>>> \Hom_\Lambda(E,Z)/t\Hom_\Lambda(E,Z) @>>> G(u) @>>> 0.
	 $$
	 Note that $F_1(u)/tF_1(u) = F(u)$ and that
	 $\Hom_\Lambda(E,Z)/t\Hom_\Lambda(E,Z)$ is a projective $\Pi_{n-1}$-module.
	 It follows that $F(u)$ coincides in $\underline \mod\ \Pi_{n-1}$
	 with $\Omega G(u).$
	\medskip
This completes the proof of Theorem 2.
	      \bigskip
	      {\bf Remark.} It should be stressed that for $n\ge 2$,
	      there cannot exist an endofunctor $\phi$
	      of $\mod \Pi_{n-1}$ such that $F
	      = \phi G$ or $\phi F= G$. For example, if we would have
	      $F = \phi G$, then the set of objects of $\Cal S(n)$ killed by
	       $G$ would be contained in the set of objects killed by $F$.
	       However $([0],[1])$ is killed by $G$, but not by $F$. Similarly,
	       $([1],[n])$ is killed by $F$, but not by $G.$

		\bigskip\bigskip
		{\bf 8. Abelian factor categories of triangulated categories}
		     \medskip
		     As a byproduct of our consideration, we see that here we deal with
		     examples of triangulated categories $\Cal T$ with an ideal $\Cal I$
		     generated by an idempotent such that $\Cal T/\Cal I$ is abelian.

		     Namely, let $\Cal T = \underline{\Cal S}(n)$ be the stable category
		     of $\Cal S(n)$, it is obtained from $\Cal S(n)$ by factoring out the ideal
		     generated by the objects $([0],[n])$ and $([n],[n])$. Now $\Cal S(n)$
		     is in a natural way (see [RS2]) a Frobenius category
		     such that $([0],[n])$ and $([n],[n])$ are the
		     only indecomposable objects which are both projective and injective, thus
		     the category $\Cal T$ is a triangulated category.

		     Let $\Cal I$ be the ideal in $\Cal T$ generated by the objects $([i],[j])$ with
		     $1\le i \le n-1$ and either $i = j$ or $j = n$. Then, as additive
		     categories, we have equivalences
		     $$
		      \Cal T/\Cal I \simeq \Cal S(n)/ \Cal U  \simeq \mod \Pi_{n-1},
		      $$
		      thus $\Cal T/\Cal I$ is an abelian category.

		      Similarly, let $\Cal J$ be the ideal in $\Cal T$
		      generated by the objects $([i],[j])$ with
		      $1\le j \le n-1$ and either $i = j$ or $i = 0$. Then, as additive
		      categories, we have equivalences
		      $$
		       \Cal T/\Cal J \simeq \Cal S(n)/ \Cal V  \simeq \mod \Pi_{n-1},
		       $$
		       thus $\Cal T/\Cal J$ is again an abelian category.

		       Here, the ideals $\Cal I, \Cal J$ of $\Cal T$ each are generated
		       by $2n-2$ indecomposable objects, whereas the rank of the Grothendieck
		       group $K_0(\Cal T/\Cal I)$ is $n-1$. One may compare this with
		       examples which involve cluster categories and cluster tilted
		       algebras, see [BMR] and [K]. Let $T$ be a cluster tilting object in
		       a cluster category $\Cal T$ say of type $Q$, where $Q$ is a directed quiver with
		       $n$ vertices. Then the factor category $\Cal T/T$ is an
		       abelian category and the rank of the Grothendieck
		       group $K_0(\Cal T/T)$ is equal to $n$, and this is also the number
		       of isomorphism classes of indecomposable direct summands of $T$.
		       
		           \bigskip\bigskip
 
{\bf 9. More objective functors $\Cal S(n) \to \mod \Pi_{n-1}$.}
     \medskip
The main part of the paper was devoted to a detailed study of two objective
functors  $\Cal S(n) \to \mod \Pi_{n-1}$, but there are also other ones.

As in the previous section, let us denote by $\Cal T$ the stable category of $\Cal S(n)$
by $\Cal T$. 
Let $\Cal I$ be again the ideal of $\Cal T$ generated by the objects
$([i],[j])$ with $1\le i \le n-1$ and either $i=j$ or $j = n$, thus 
there is an equivalence $\zeta\:\Cal T/\Cal I \to  \mod \Pi_{n-1}$.

Thus, {\it let $\eta$ be any autoequivalence of $\Cal T$. Then the composition 
$$
  \Cal S(n) @>>> \Cal T @>\eta>>  \Cal T @>>> \Cal T/\Cal I 
  @>\zeta>> \mod \Pi_{n-1}
$$ 
(where the first and the third functor are the canonical projections) clearly is
a full, dense, objective functor.} 
	\medskip
We recall from [RS1] that the stable category $\Cal T$ has non-trivial 
autoequivalences, for example the endofunctor induced by the Auslander-Reiten
translation $\tau$ is an autoequivalence of order $6$.

	\bigskip\bigskip 

{\bf 10. Auslander-Reiten orbits.}
     \medskip
     {\bf Some Auslander-Reiten orbits in the category $\Cal S(n).$}
     Let $\tau$ be the Auslander-Reiten translation in $\Cal S(n)$.
     The paper [RS1] describes in detail how to obtain for the pair
     $(X,Y)$ in $\Cal S(n)$ the pair
     $\tau(X,Y)$.

     We are interested in some of the objects of the form $([i],[j])$ (with $0\le i \le j \le n$).
     The objects $([n],[n])$ and $([0],[n])$ are projective-injective, thus they are
     sent to zero by $\tau.$
     The following assertions for $1\le i \le n\!-\!1$
     $$
     \align
      \tau([0],[i]) &= ([i],[i])\cr
       \tau([i],[i]) &= ([i],[n])\cr
        \tau([i],[n]) &= ([0],[n\!-\!i])\cr
	\endalign
	$$
	are easily verified. Of course, $1\le i \le n\!-\!1$ implies that also
	$1\le n\!-\!i \le n\!-\!1$,
	thus we see that the set of objects
	of the form $([0],[i]),([i],[i]),([i],[n])$ with $1\le i\le n\!-\!1$ is closed under $\tau$.
	Let us present the corresponding
	parts of Auslander-Reiten components. First of all, for $2 \le i < \frac n2$, we deal with a
	$\tau$-orbit of length 6:
	$$
	
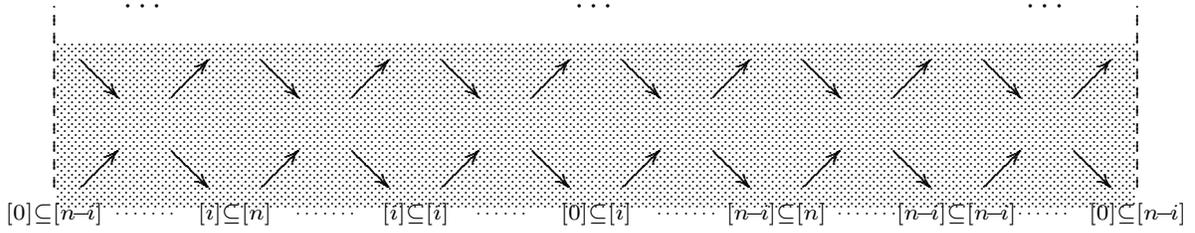
\hbox{\beginpicture
	\setcoordinatesystem units <1cm,1cm>
	\put{} at 0 .7
	\multiput{} at -1 0  13 0 /

	\arr{0.3 0.3}{0.7 0.7}
	\arr{2.3 0.3}{2.7 0.7}
	\arr{4.3 0.3}{4.7 0.7}
	\arr{1.3 0.7}{1.7 0.3}
	\arr{3.3 0.7}{3.7 0.3}
	\arr{5.3 0.7}{5.7 0.3}

	\setdots <1mm>
	\plot .7 0  1.3 0 /
	\plot 2.7 0  3.3 0 /
	\plot 4.7 0  5.2 0 /

	\setsolid
	\arr{6.3 0.3}{6.7 0.7}
	\arr{8.3 0.3}{8.7 0.7}
	\arr{10.3 0.3}{10.7 0.7}
	\arr{7.3 0.7}{7.7 0.3}
	\arr{9.3 0.7}{9.7 0.3}
	\arr{11.3 0.7}{11.7 0.3}

	\setdots <1mm>
	\plot 6.7 0  7.3 0 /
	\plot 8.7 0  9.3 0 /
	\plot 10.7 0  11.2 0 /

	\put{$\ssize[i]\subseteq[i]$} at 4 0
	\put{$\ssize[i]\subseteq[n]$} at 2 0
	\multiput{$\ssize[0]\subseteq[n\!-\!i]$} at 0 0  12 0 /

	\put{$\ssize[n\!-\!i]\subseteq[n\!-\!i]$} at 10 0
	\put{$\ssize[n\!-\!i]\subseteq[n]$} at 8 0
	\put{$\ssize[0]\subseteq[i]$} at 6 0

	\put{} at 0 2.3

	\setshadegrid span <.4mm>
	\vshade  0 0.05 1.9 <,,,>  12 0.05 1.9  /
	\setsolid
	\arr{0.3 1.7}{0.7 1.3}
	\arr{1.3 1.3}{1.7 1.7}
	\arr{2.3 1.7}{2.7 1.3}
	\arr{3.3 1.3}{3.7 1.7}
	\arr{4.3 1.7}{4.7 1.3}
	\arr{5.3 1.3}{5.7 1.7}
	\arr{6.3 1.7}{6.7 1.3}
	\arr{7.3 1.3}{7.7 1.7}
	\arr{8.3 1.7}{8.7 1.3}
	\arr{9.3 1.3}{9.7 1.7}
	\arr{10.3 1.7}{10.7 1.3}
	\arr{11.3 1.3}{11.7 1.7}
	\setdashes <1mm>
	\plot 0 0.3  0 2.3 /
	\plot 12 0.3  12 2.3 /
	\multiput{$\cdots$} at 1 2.3  6 2.3  11 2.3 /

	\endpicture}
	$$
	If $n\ge 4$ is even and $i =\frac n2$, there is the following orbit of length 3:
	$$
	
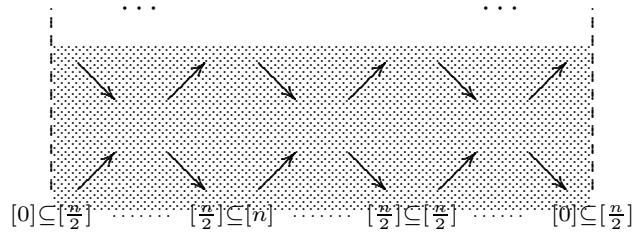
\hbox{\beginpicture
	\setcoordinatesystem units <1cm,1cm>
	\put{} at 0 .7

	\arr{0.3 0.3}{0.7 0.7}
	\arr{2.3 0.3}{2.7 0.7}
	\arr{4.3 0.3}{4.7 0.7}
	\arr{1.3 0.7}{1.7 0.3}
	\arr{3.3 0.7}{3.7 0.3}
	\arr{5.3 0.7}{5.7 0.3}

	\setdots <1mm>
	\plot .7 0  1.3 0 /
	\plot 2.7 0  3.3 0 /
	\plot 4.7 0  5.2 0 /

	\put{$\ssize[\frac n2]\subseteq[\frac n2]$} at 4 0
	\put{$\ssize[\frac n2]\subseteq[n]$} at 2 0
	\multiput{$\ssize[0]\subseteq[\frac n2]$} at 0 0  6 0 /

	\put{} at 0 2.3

	\setshadegrid span <.4mm>
	\vshade  0 0.05 1.9 <,,,>  6 0.05 1.9  /
	\setsolid
	\arr{0.3 1.7}{0.7 1.3}
	\arr{1.3 1.3}{1.7 1.7}
	\arr{2.3 1.7}{2.7 1.3}
	\arr{3.3 1.3}{3.7 1.7}
	\arr{4.3 1.7}{4.7 1.3}
	\arr{5.3 1.3}{5.7 1.7}
	\setdashes <1mm>
	\plot 0 0.3  0 2.3 /
	\plot 6 0.3  6 2.3 /
	\multiput{$\cdots$} at 1 2.3  5 2.3  /

	\endpicture}
	$$
	Finally, for $i=1$ we get:
	$$
	
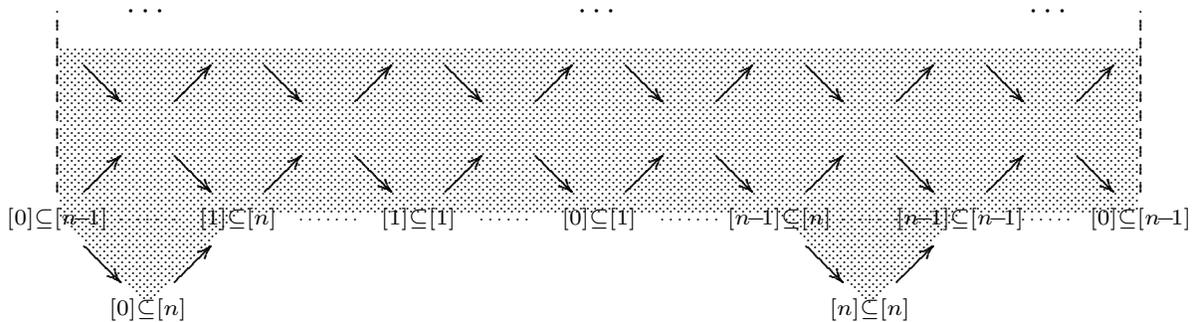
\hbox{\beginpicture
	\setcoordinatesystem units <1cm,1cm>
	\put{} at 0 .7
	\multiput{} at -1 0  13 0 /

	\put{$\ssize[1]\subseteq[1]$} at 4 0
	\put{$\ssize[1]\subseteq[n]$} at 2 0
	\multiput{$\ssize[0]\subseteq[n\!-\!1]$} at 0 0  12 0 /
	\put{$\ssize[0]\subseteq[n]$} at 1 -1
	\arr{0.3 0.3}{0.7 0.7}
	\arr{2.3 0.3}{2.7 0.7}
	\arr{4.3 0.3}{4.7 0.7}
	\arr{1.3 0.7}{1.7 0.3}
	\arr{3.3 0.7}{3.7 0.3}
	\arr{5.3 0.7}{5.7 0.3}

	\arr{0.3 -.3}{0.7 -.7}
	\arr{1.3 -.7}{1.7 -.3}

	\setdots <1mm>

	\plot .7 0  1.3 0 /
	\plot 2.7 0  3.3 0 /
	\plot 4.7 0  5.2 0 /

	\setsolid
	\put{$\ssize[n\!-\!1]\subseteq[n\!-\!1]$} at 10 0
	\put{$\ssize[n\!-\!1]\subseteq[n]$} at 8 0
	\put{$\ssize[0]\subseteq[1]$} at 6 0
	\put{$\ssize[n]\subseteq[n]$} at 9 -1
	\arr{6.3 0.3}{6.7 0.7}
	\arr{8.3 0.3}{8.7 0.7}
	\arr{10.3 0.3}{10.7 0.7}
	\arr{7.3 0.7}{7.7 0.3}
	\arr{9.3 0.7}{9.7 0.3}
	\arr{11.3 0.7}{11.7 0.3}

	\arr{8.3 -.3}{8.7 -.7}
	\arr{9.3 -.7}{9.7 -.3}
	\setdots <1mm>
	\plot 6.7 0  7.3 0 /
	\plot 8.7 0  9.3 0 /
	\plot 10.7 0  11.2 0 /

	\put{} at 0 2.3

	\setshadegrid span <.4mm>
	\vshade  0 0.05 1.9 <,z,,> 1 -1 1.9 <z,z,,>  2 0.05 1.9
	  <z,z,,> 8 0.05 1.9 <z,z,,> 9 -1 1.9
	   <z,z,,> 10 0.05 1.9  <z,,,>  12 0.05 1.9  /
	   \setsolid
	   \arr{0.3 1.7}{0.7 1.3}
	   \arr{1.3 1.3}{1.7 1.7}
	   \arr{2.3 1.7}{2.7 1.3}
	   \arr{3.3 1.3}{3.7 1.7}
	   \arr{4.3 1.7}{4.7 1.3}
	   \arr{5.3 1.3}{5.7 1.7}
	   \arr{6.3 1.7}{6.7 1.3}
	   \arr{7.3 1.3}{7.7 1.7}
	   \arr{8.3 1.7}{8.7 1.3}
	   \arr{9.3 1.3}{9.7 1.7}
	   \arr{10.3 1.7}{10.7 1.3}
	   \arr{11.3 1.3}{11.7 1.7}
	   \setdashes <1mm>
	   \plot 0 0.3  0 2.3 /
	   \plot 12 0.3  12 2.3 /
	   \multiput{$\cdots$} at 1 2.3  6 2.3  11 2.3 /

	   \endpicture}
	   $$
		\bigskip\bigskip
		\noindent
		{\bf The corresponding $\tau$-orbits in the category $\Cal F(n).$}
		First, those for $2\le i < \frac n2$:
		$$
		
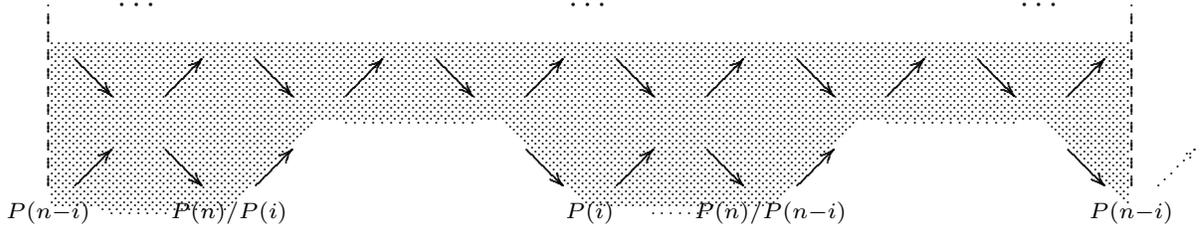
\hbox{\beginpicture
		\setcoordinatesystem units <1cm,1cm>
		\put{} at 0 .7
		\multiput{} at -1 0  13 0 /

		\arr{0.3 0.3}{0.7 0.7}
		\arr{2.3 0.3}{2.7 0.7}
		\arr{1.3 0.7}{1.7 0.3}
		\arr{5.3 0.7}{5.7 0.3}

		\setdots <1mm>

		\plot .7 0  1.3 0 /

		\setsolid
		\arr{6.3 0.3}{6.7 0.7}
		\arr{8.3 0.3}{8.7 0.7}
		\arr{7.3 0.7}{7.7 0.3}
		\arr{11.3 0.7}{11.7 0.3}

		\setdots <1mm>
		\plot 6.7 0  7.3 0 /

		\put{$\ssize P(n)/P(i)$} at 2 0
		\multiput{$\ssize P(n-i)$} at 0 0  12 0 /

		\put{$\ssize P(n)/P(n-i)$} at 8 0
		\put{$\ssize P(i)$} at 6 0
		\put{} at -1 0
		\arr{12.3 0.3}{12.7 0.7}

		\plot 3.3 1  4.7 1 /
		\plot 9.3 1  10.7 1 /
		\put{} at 0 2.3

		\setshadegrid span <.4mm>
		\vshade  0 0.05 1.9 <,z,,> 2 0 1.9 <z,z,,>  3  1 1.9
		 <z,z,,>  5 1 1.9 <z,z,,> 6 0.05 1.9 <z,z,,> 8 0.05 1.9 <z,z,,> 9 1 1.9
		  <z,z,,> 11 1 1.9  <z,,,>  12 0.05 1.9  /
		  \setsolid
		  \arr{0.3 1.7}{0.7 1.3}
		  \arr{1.3 1.3}{1.7 1.7}
		  \arr{2.3 1.7}{2.7 1.3}
		  \arr{3.3 1.3}{3.7 1.7}
		  \arr{4.3 1.7}{4.7 1.3}
		  \arr{5.3 1.3}{5.7 1.7}
		  \arr{6.3 1.7}{6.7 1.3}
		  \arr{7.3 1.3}{7.7 1.7}
		  \arr{8.3 1.7}{8.7 1.3}
		  \arr{9.3 1.3}{9.7 1.7}
		  \arr{10.3 1.7}{10.7 1.3}
		  \arr{11.3 1.3}{11.7 1.7}
		  \setdashes <1mm>
		  \plot 0 0.3  0 2.3 /
		  \plot 12 0.3  12 2.3 /
		  \multiput{$\cdots$} at 1 2.3  6 2.3  11 2.3 /

		  \endpicture}
		  $$
		  Second, for $n\ge 4$ even and $i=\frac n2$:
		  $$
		  
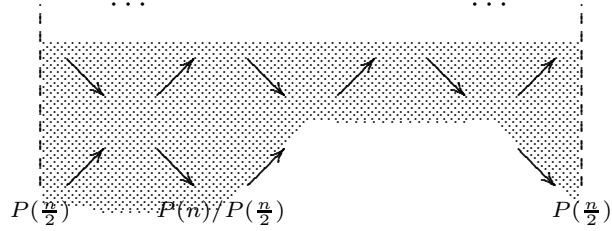
\hbox{\beginpicture
		  \setcoordinatesystem units <1cm,1cm>
		  \put{} at 0 .7
		  \arr{0.3 0.3}{0.7 0.7}
		  \arr{2.3 0.3}{2.7 0.7}
		  \arr{1.3 0.7}{1.7 0.3}
		  \arr{5.3 0.7}{5.7 0.3}

		  \setdots <1mm>
		  \plot .7 0  1.3 0 /

		  \put{$\ssize P(n)/P(\frac n2)$} at 2 0
		  \multiput{$\ssize P(\frac n2)$} at 0 0  6 0 /

		  \plot 3.3 1  4.7 1 /
		  \put{} at 0 2.3

		  \setshadegrid span <.4mm>
		  \vshade  0 0.05 1.9 <,z,,> 2 0 1.9 <z,z,,>  3  1 1.9
		   <z,z,,>  5 1 1.9 <z,,,> 6 0.05 1.9   /
		   \setsolid
		   \arr{0.3 1.7}{0.7 1.3}
		   \arr{1.3 1.3}{1.7 1.7}
		   \arr{2.3 1.7}{2.7 1.3}
		   \arr{3.3 1.3}{3.7 1.7}
		   \arr{4.3 1.7}{4.7 1.3}
		   \arr{5.3 1.3}{5.7 1.7}
		   \setdashes <1mm>
		   \plot 0 0.3  0 2.3 /
		   \plot 6 0.3  6 2.3 /
		   \multiput{$\cdots$} at 1 2.3  5 2.3 /

		   \endpicture}
		   $$
		   And finally, for $i=1$ we get:
		   $$
		   
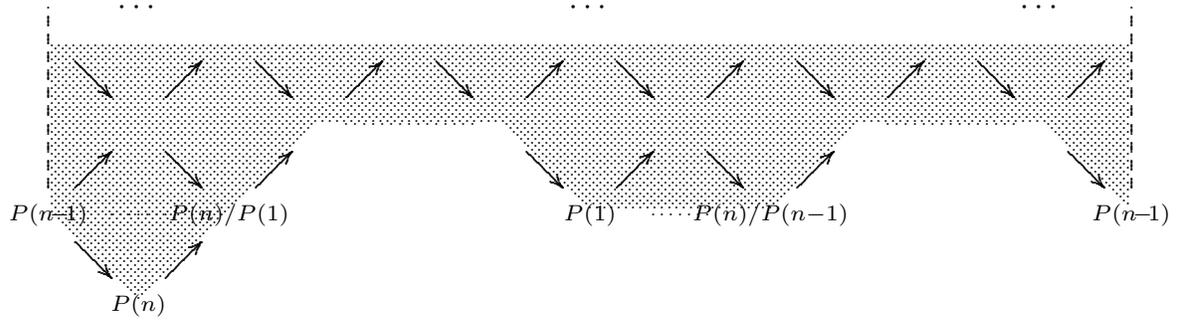
\hbox{\beginpicture
		   \setcoordinatesystem units <1cm,1cm>
		   \put{} at 0 .7
		   \multiput{} at -1 0  13 0 /

		   \put{$\ssize P(n)/P(1)$} at 2 0
		   \multiput{$\ssize P(n\!-\!1)$} at 0 0  12 0 /
		   \put{$\ssize P(n)$} at 1 -1
		   \arr{0.3 0.3}{0.7 0.7}
		   \arr{2.3 0.3}{2.7 0.7}
		   \arr{1.3 0.7}{1.7 0.3}
		   \arr{5.3 0.7}{5.7 0.3}

		   \arr{0.3 -.3}{0.7 -.7}
		   \arr{1.3 -.7}{1.7 -.3}

		   \setdots <1mm>
		   \plot .7 0  1.3 0 /

		   \setsolid
		   \put{$\ssize P(n)/P(n-1)$} at 8 0
		   \put{$\ssize P(1)$} at 6 0
		   \arr{6.3 0.3}{6.7 0.7}
		   \arr{8.3 0.3}{8.7 0.7}
		   \arr{7.3 0.7}{7.7 0.3}
		   \arr{11.3 0.7}{11.7 0.3}

		   \setdots <1mm>
		   \plot 6.7 0  7.3 0 /
		   \plot 3.3 1  4.7 1 /
		   \plot 9.3 1  10.7 1 /
		   \put{} at 0 2.3
		   \setshadegrid span <.4mm>
		   \vshade  0 0.05 1.9 <,z,,> 1 -1 1.9 <z,z,,>  3  1 1.9
		    <z,z,,>  5 1 1.9 <z,z,,> 6 0.05 1.9 <z,z,,> 8 0.05 1.9 <z,z,,> 9 1 1.9
		     <z,z,,> 11 1 1.9  <z,,,>  12 0.05 1.9  /
		     \setsolid
		     \arr{0.3 1.7}{0.7 1.3}
		     \arr{1.3 1.3}{1.7 1.7}
		     \arr{2.3 1.7}{2.7 1.3}
		     \arr{3.3 1.3}{3.7 1.7}
		     \arr{4.3 1.7}{4.7 1.3}
		     \arr{5.3 1.3}{5.7 1.7}
		     \arr{6.3 1.7}{6.7 1.3}
		     \arr{7.3 1.3}{7.7 1.7}
		     \arr{8.3 1.7}{8.7 1.3}
		     \arr{9.3 1.3}{9.7 1.7}
		     \arr{10.3 1.7}{10.7 1.3}
		     \arr{11.3 1.3}{11.7 1.7}
		     \setdashes <1mm>
		     \plot 0 0.3  0 2.3 /
		     \plot 12 0.3  12 2.3 /
		     \multiput{$\cdots$} at 1 2.3  6 2.3  11 2.3 /
		     \endpicture}
		     $$
			\bigskip
			{\bf Remark.} The paper [RS2] describes in detail the Auslander-Reiten
			quivers of the categories $\Cal S(n)$ with $1\le n \le 6$. Similarly,
			in [DR] the Auslander-Reiten quivers of the categories $\Cal F(n)$
			with $2\le n \le 5$ are presented. As the functor $F_1$ shows, the Auslander-Reiten
			quiver of $\Cal F(n)$ can be obtained from that of $\Cal S(n)$ by just
			deleting some vertices, thus it is easy to obtain the illustrations
			presented in [DR]
			from those in [RS2].
			The deletion process explains also some of the
			features of the shape of the Auslander-Reiten quiver of $\Cal F(n)$:
			Of course, there is precisely one projective-injective vertex, namely
			the module $P(n) = I(n) = T(1)$. The remaining transjective orbits contain precisely
			two vertices, namely $T(n+1-i)$ and $P(i) = \tau_{\Cal F(n)} T(n+1-i)$,
			here $1\le i\le n-1$. As we now know, this concerns the $\tau$-orbit of
			$\Cal S(n)$ which contains the pairs $([i],[i])$ and $([n-i],[n-i])$:
			both pairs are killed by the functor $F_1$, but in-between the
			Auslander-Reiten sequence starting with $([0],[i])$ and ending in
			$([n-i],[n])$ is not touched and it yields under $F_1$ the Auslander-Reiten sequence
			starting with $P(i)$ and ending in $T(n+1-i).$
				  \bigskip\bigskip
				  \noindent
				  {\bf The corresponding $\tau$-orbits
				  of the category $\Cal S(n)/\Cal U \simeq \mod \Pi_{n-1}$.}

				  \noindent
				  For $2\le i < \frac n2$:
				  $$
				  
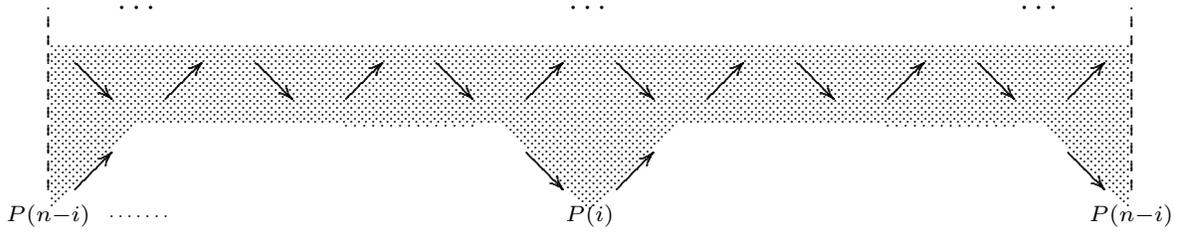
\hbox{\beginpicture
				  \setcoordinatesystem units <1cm,1cm>
				  \put{} at 0 .7
				  \multiput{} at -1 0  13 0 /

				  \arr{0.3 0.3}{0.7 0.7}
				  \arr{5.3 0.7}{5.7 0.3}

				  \setdots <1mm>
				  \plot .7 0  1.3 0 /

				  \setsolid
				  \arr{6.3 0.3}{6.7 0.7}
				  \arr{11.3 0.7}{11.7 0.3}

				  \multiput{$\ssize P(n-i)$} at 0 0  12 0 /

				  \put{$\ssize P(i)$} at 6 0
				  \put{} at -1 0

				  \setdots <1mm>
				  \plot 3.3 1  4.7 1 /
				  \plot 9.3 1  10.7 1 /
				  \put{} at 0 2.3

				  \setshadegrid span <.4mm>
				  \vshade  0 0.05 1.9 <,z,,> 1 1 1.9
				   <z,z,,>  5 1 1.9 <z,z,,> 6 0 1.9 <z,z,,> 7 1 1.9
				    <z,z,,> 11 1 1.9  <z,,,>  12 0.05 1.9  /
				    \setsolid
				    \arr{0.3 1.7}{0.7 1.3}
				    \arr{1.3 1.3}{1.7 1.7}
				    \arr{2.3 1.7}{2.7 1.3}
				    \arr{3.3 1.3}{3.7 1.7}
				    \arr{4.3 1.7}{4.7 1.3}
				    \arr{5.3 1.3}{5.7 1.7}
				    \arr{6.3 1.7}{6.7 1.3}
				    \arr{7.3 1.3}{7.7 1.7}
				    \arr{8.3 1.7}{8.7 1.3}
				    \arr{9.3 1.3}{9.7 1.7}
				    \arr{10.3 1.7}{10.7 1.3}
				    \arr{11.3 1.3}{11.7 1.7}
				    \setdashes <1mm>
				    \plot 0 0.3  0 2.3 /
				    \plot 12 0.3  12 2.3 /
				    \multiput{$\cdots$} at 1 2.3  6 2.3  11 2.3 /

				    \endpicture}
				    $$
				    For $n\ge 4$ even and $i=\frac n2$:
				    $$
				    
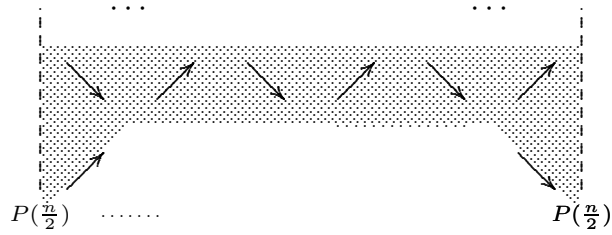
\hbox{\beginpicture
				    \setcoordinatesystem units <1cm,1cm>
				    \put{} at 0 .7

				    \arr{0.3 0.3}{0.7 0.7}
				    \arr{5.3 0.7}{5.7 0.3}

				    \setdots <1mm>
				    \plot .7 0  1.3 0 /

				    \multiput{$\ssize P(\frac n2)$} at 0 0  6 0 /

				    \put{$\ssize P(\frac n2)$} at 6 0

				    \plot 3.3 1  4.7 1 /
				    \put{} at 0 2.3

				    \setshadegrid span <.4mm>
				    \vshade  0 0.05 1.9 <,z,,> 1 1 1.9
				     <z,z,,>  5 1 1.9 <z,,,> 6 0 1.9  /
				     \setsolid
				     \arr{0.3 1.7}{0.7 1.3}
				     \arr{1.3 1.3}{1.7 1.7}
				     \arr{2.3 1.7}{2.7 1.3}
				     \arr{3.3 1.3}{3.7 1.7}
				     \arr{4.3 1.7}{4.7 1.3}
				     \arr{5.3 1.3}{5.7 1.7}
				     \setdashes <1mm>
				     \plot 0 0.3  0 2.3 /
				     \plot 6 0.3  6 2.3 /
				     \multiput{$\cdots$} at 1 2.3  5 2.3 /

				     \endpicture}
				     $$
				     And finally, for $i=1$ we get:
				     $$
				     
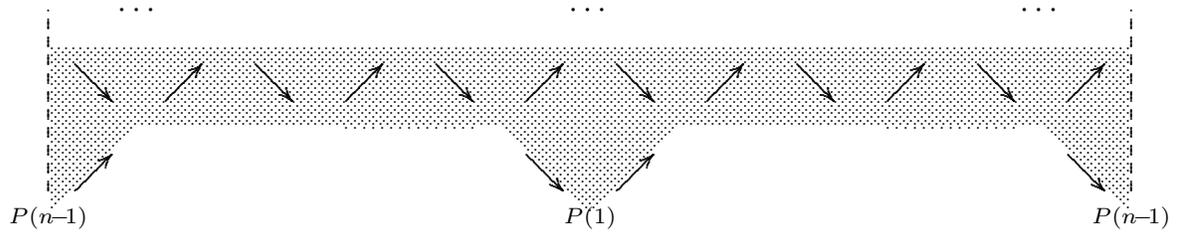
\hbox{\beginpicture
				     \setcoordinatesystem units <1cm,1cm>
				     \put{} at 0 .7
				     \multiput{} at -1 0  13 0 /

				     \multiput{$\ssize P(n\!-\!1)$} at 0 0  12 0 /

				     \arr{0.3 0.3}{0.7 0.7}
				     \arr{5.3 0.7}{5.7 0.3}

				     \setdots <1mm>
				     \setsolid

				     \put{$\ssize P(1)$} at 6 0

				     \arr{6.3 0.3}{6.7 0.7}
				     \arr{11.3 0.7}{11.7 0.3}

				     \setdots <1mm>

				     \plot 3.3 1  4.7 1 /
				     \plot 9.3 1  10.7 1 /
				     \put{} at 0 2.3
				     \setshadegrid span <.4mm>
				     \vshade  0 0.05 1.9 <,z,,> 1 1 1.9
				      <z,z,,>  5 1 1.9 <z,z,,> 6 0 1.9 <z,z,,> 7 1 1.9
				       <z,z,,> 11 1 1.9  <z,,,>  12 0.05 1.9  /

				       \setsolid
				       \arr{0.3 1.7}{0.7 1.3}
				       \arr{1.3 1.3}{1.7 1.7}
				       \arr{2.3 1.7}{2.7 1.3}
				       \arr{3.3 1.3}{3.7 1.7}
				       \arr{4.3 1.7}{4.7 1.3}
				       \arr{5.3 1.3}{5.7 1.7}
				       \arr{6.3 1.7}{6.7 1.3}
				       \arr{7.3 1.3}{7.7 1.7}
				       \arr{8.3 1.7}{8.7 1.3}
				       \arr{9.3 1.3}{9.7 1.7}
				       \arr{10.3 1.7}{10.7 1.3}
				       \arr{11.3 1.3}{11.7 1.7}
				       \setdashes <1mm>
				       \plot 0 0.3  0 2.3 /
				       \plot 12 0.3  12 2.3 /
				       \multiput{$\cdots$} at 1 2.3  6 2.3  11 2.3 /
				       \endpicture}
				       $$
					\bigskip\bigskip
					{\bf A simultaneous view.} Let us draw again the relevant components of
					$\Cal S(n)$ and use the shading in order to illustrate
					what remains when we delete the objects in $\Cal U$.
					Note that under the functor $F$, we have
					$$
					\alignat 2
					 F([i],[i]) &= 0     &&1\le i \le n, \cr
					  F([i],[n]) &= 0 &\qquad\text{for}\qquad &0\le i \le n\!-\!1, \cr
					   F([0],[j]) &= P_{\Pi}(j)  &&1\le j\le n\!-\!1. \cr
					   \endalignat
					   $$

						\noindent
						First of all, for $2 \le i < \frac n2$, two objects of the $\tau_{\Cal S}$-orbit
						of $([i],[i])$ survive:
						$$
						
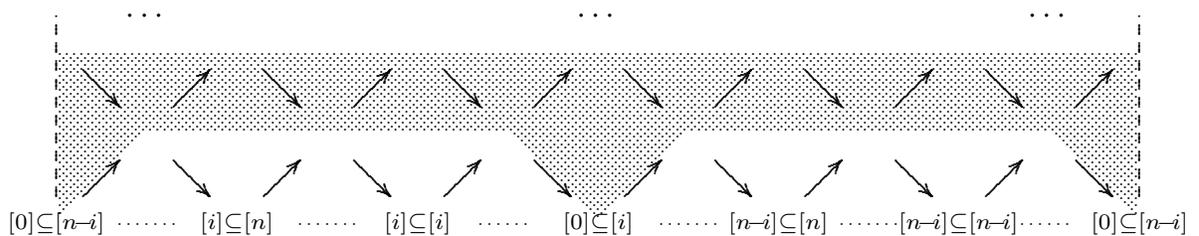
\hbox{\beginpicture
						\setcoordinatesystem units <1cm,1cm>
						\put{} at 0 .7
						\multiput{} at -1 0  13 0 /

						\arr{0.3 0.3}{0.7 0.7}
						\arr{2.3 0.3}{2.7 0.7}
						\arr{4.3 0.3}{4.7 0.7}
						\arr{1.3 0.7}{1.7 0.3}
						\arr{3.3 0.7}{3.7 0.3}
						\arr{5.3 0.7}{5.7 0.3}

						\setdots <1mm>
						\plot .7 0  1.3 0 /
						\plot 2.7 0  3.3 0 /
						\plot 4.7 0  5.2 0 /

						\setsolid
						\arr{6.3 0.3}{6.7 0.7}
						\arr{8.3 0.3}{8.7 0.7}
						\arr{10.3 0.3}{10.7 0.7}
						\arr{7.3 0.7}{7.7 0.3}
						\arr{9.3 0.7}{9.7 0.3}
						\arr{11.3 0.7}{11.7 0.3}

						\setdots <1mm>
						\plot 6.7 0  7.3 0 /
						\plot 8.7 0  9.3 0 /
						\plot 10.7 0  11.2 0 /

						\put{$\ssize[i]\subseteq[i]$} at 4 0
						\put{$\ssize[i]\subseteq[n]$} at 2 0
						\multiput{$\ssize[0]\subseteq[n\!-\!i]$} at 0 0  12 0 /

						\put{$\ssize[n\!-\!i]\subseteq[n\!-\!i]$} at 10 0
						\put{$\ssize[n\!-\!i]\subseteq[n]$} at 8 0
						\put{$\ssize[0]\subseteq[i]$} at 6 0

						\put{} at 0 2.3

						\setshadegrid span <.4mm>
						\vshade  0 0.05 1.9 <,z,,> 1 1 1.9
						 <z,z,,>  5 1 1.9 <z,z,,> 6 0 1.9 <z,z,,> 7 1 1.9
						  <z,z,,> 11 1 1.9  <z,,,>  12 0.05 1.9  /
						  \setsolid
						  \arr{0.3 1.7}{0.7 1.3}
						  \arr{1.3 1.3}{1.7 1.7}
						  \arr{2.3 1.7}{2.7 1.3}
						  \arr{3.3 1.3}{3.7 1.7}
						  \arr{4.3 1.7}{4.7 1.3}
						  \arr{5.3 1.3}{5.7 1.7}
						  \arr{6.3 1.7}{6.7 1.3}
						  \arr{7.3 1.3}{7.7 1.7}
						  \arr{8.3 1.7}{8.7 1.3}
						  \arr{9.3 1.3}{9.7 1.7}
						  \arr{10.3 1.7}{10.7 1.3}
						  \arr{11.3 1.3}{11.7 1.7}
						  \setdashes <1mm>
						  \plot 0 0.3  0 2.3 /
						  \plot 12 0.3  12 2.3 /
						  \multiput{$\cdots$} at 1 2.3  6 2.3  11 2.3 /

						  \endpicture}
						  $$
							\bigskip
							\noindent
							If $n\ge 4$ is even and $i =\frac n2$, the pair $([0],[i])$ is the
							only object in the $\tau_{\Cal S}$-orbit of $([i],[i])$ which survives:
							$$
							
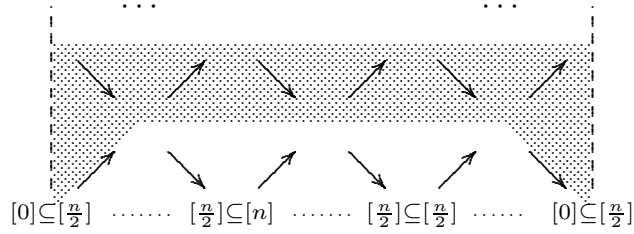
\hbox{\beginpicture
							\setcoordinatesystem units <1cm,1cm>
							\put{} at 0 .7

							\arr{0.3 0.3}{0.7 0.7}
							\arr{2.3 0.3}{2.7 0.7}
							\arr{4.3 0.3}{4.7 0.7}
							\arr{1.3 0.7}{1.7 0.3}
							\arr{3.3 0.7}{3.7 0.3}
							\arr{5.3 0.7}{5.7 0.3}

							\setdots <1mm>
							\plot .7 0  1.3 0 /
							\plot 2.7 0  3.3 0 /
							\plot 4.7 0  5.2 0 /

							\put{$\ssize[\frac n2]\subseteq[\frac n2]$} at 4 0
							\put{$\ssize[\frac n2]\subseteq[n]$} at 2 0
							\multiput{$\ssize[0]\subseteq[\frac n2]$} at 0 0  6 0 /

							\put{} at 0 2.3

							\setshadegrid span <.4mm>
							\vshade  0 0.05 1.9 <,z,,> 1 1 1.9
							 <z,z,,>  5 1 1.9 <z,,,> 6 0 1.9  /

							 \setsolid
							 \arr{0.3 1.7}{0.7 1.3}
							 \arr{1.3 1.3}{1.7 1.7}
							 \arr{2.3 1.7}{2.7 1.3}
							 \arr{3.3 1.3}{3.7 1.7}
							 \arr{4.3 1.7}{4.7 1.3}
							 \arr{5.3 1.3}{5.7 1.7}
							 \setdashes <1mm>
							 \plot 0 0.3  0 2.3 /
							 \plot 6 0.3  6 2.3 /
							 \multiput{$\cdots$} at 1 2.3  5 2.3  /

							 \endpicture}
							 $$
								\bigskip
								\noindent
								Finally, for $i=1$, again two objects in the $\tau_{\Cal S}$-orbit of
								$([1],[1])$ survive:
								$$
								
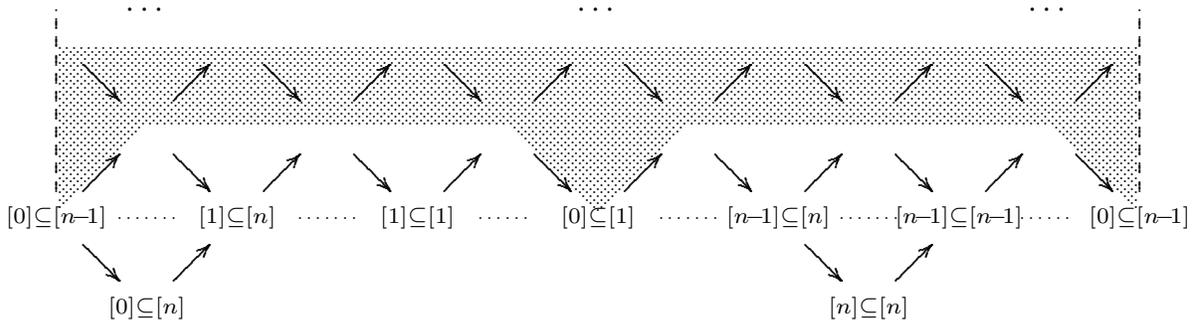
\hbox{\beginpicture
								\setcoordinatesystem units <1cm,1cm>
								\put{} at 0 .7
								\multiput{} at -1 0  13 0 /

								\put{$\ssize[1]\subseteq[1]$} at 4 0
								\put{$\ssize[1]\subseteq[n]$} at 2 0
								\multiput{$\ssize[0]\subseteq[n\!-\!1]$} at 0 0  12 0 /
								\put{$\ssize[0]\subseteq[n]$} at 1 -1
								\arr{0.3 0.3}{0.7 0.7}
								\arr{2.3 0.3}{2.7 0.7}
								\arr{4.3 0.3}{4.7 0.7}
								\arr{1.3 0.7}{1.7 0.3}
								\arr{3.3 0.7}{3.7 0.3}
								\arr{5.3 0.7}{5.7 0.3}

								\arr{0.3 -.3}{0.7 -.7}
								\arr{1.3 -.7}{1.7 -.3}

								\setdots <1mm>

								\plot .7 0  1.3 0 /
								\plot 2.7 0  3.3 0 /
								\plot 4.7 0  5.2 0 /

								\setsolid
								\put{$\ssize[n\!-\!1]\subseteq[n\!-\!1]$} at 10 0
								\put{$\ssize[n\!-\!1]\subseteq[n]$} at 8 0
								\put{$\ssize[0]\subseteq[1]$} at 6 0
								\put{$\ssize[n]\subseteq[n]$} at 9 -1
								\arr{6.3 0.3}{6.7 0.7}
								\arr{8.3 0.3}{8.7 0.7}
								\arr{10.3 0.3}{10.7 0.7}
								\arr{7.3 0.7}{7.7 0.3}
								\arr{9.3 0.7}{9.7 0.3}
								\arr{11.3 0.7}{11.7 0.3}

								\arr{8.3 -.3}{8.7 -.7}
								\arr{9.3 -.7}{9.7 -.3}
								\setdots <1mm>
								\plot 6.7 0  7.3 0 /
								\plot 8.7 0  9.3 0 /
								\plot 10.7 0  11.2 0 /

								\put{} at 0 2.3

								\setshadegrid span <.4mm>
								\vshade  0 0.05 1.9 <,z,,> 1 1 1.9
								 <z,z,,>  5 1 1.9 <z,z,,> 6 0 1.9 <z,z,,> 7 1 1.9
								  <z,z,,> 11 1 1.9  <z,,,>  12 0.05 1.9  /
								  \setsolid
								  \arr{0.3 1.7}{0.7 1.3}
								  \arr{1.3 1.3}{1.7 1.7}
								  \arr{2.3 1.7}{2.7 1.3}
								  \arr{3.3 1.3}{3.7 1.7}
								  \arr{4.3 1.7}{4.7 1.3}
								  \arr{5.3 1.3}{5.7 1.7}
								  \arr{6.3 1.7}{6.7 1.3}
								  \arr{7.3 1.3}{7.7 1.7}
								  \arr{8.3 1.7}{8.7 1.3}
								  \arr{9.3 1.3}{9.7 1.7}
								  \arr{10.3 1.7}{10.7 1.3}
								  \arr{11.3 1.3}{11.7 1.7}
								  \setdashes <1mm>
								  \plot 0 0.3  0 2.3 /
								  \plot 12 0.3  12 2.3 /
								  \multiput{$\cdots$} at 1 2.3  6 2.3  11 2.3 /

								  \endpicture}
								  $$
									\bigskip\bigskip
									{\bf The relevant components of $\Cal S(n)/\Cal V$}.
									In the same way as we have presented components of $\Cal S(n)$ shading the parts
									which remain after deleting $\Cal U$, we now
									show what remains from these components when we remove $\Cal V$.

									First of all, for $2 \le i < \frac n2$, two objects in the $\tau_{\Cal S}$-orbit
									of $([i],[i])$ survive:
									$$
									
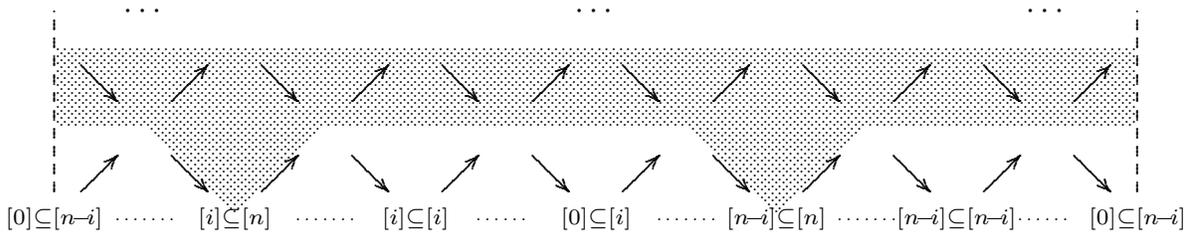
\hbox{\beginpicture
									\setcoordinatesystem units <1cm,1cm>
									\put{} at 0 .7
									\multiput{} at -1 0  13 0 /

									\arr{0.3 0.3}{0.7 0.7}
									\arr{2.3 0.3}{2.7 0.7}
									\arr{4.3 0.3}{4.7 0.7}
									\arr{1.3 0.7}{1.7 0.3}
									\arr{3.3 0.7}{3.7 0.3}
									\arr{5.3 0.7}{5.7 0.3}

									\setdots <1mm>
									\plot .7 0  1.3 0 /
									\plot 2.7 0  3.3 0 /
									\plot 4.7 0  5.2 0 /

									\setsolid
									\arr{6.3 0.3}{6.7 0.7}
									\arr{8.3 0.3}{8.7 0.7}
									\arr{10.3 0.3}{10.7 0.7}
									\arr{7.3 0.7}{7.7 0.3}
									\arr{9.3 0.7}{9.7 0.3}
									\arr{11.3 0.7}{11.7 0.3}

									\setdots <1mm>
									\plot 6.7 0  7.3 0 /
									\plot 8.7 0  9.3 0 /
									\plot 10.7 0  11.2 0 /

									\put{$\ssize[i]\subseteq[i]$} at 4 0
									\put{$\ssize[i]\subseteq[n]$} at 2 0
									\multiput{$\ssize[0]\subseteq[n\!-\!i]$} at 0 0  12 0 /

									\put{$\ssize[n\!-\!i]\subseteq[n\!-\!i]$} at 10 0
									\put{$\ssize[n\!-\!i]\subseteq[n]$} at 8 0
									\put{$\ssize[0]\subseteq[i]$} at 6 0

									\put{} at 0 2.3

									\setshadegrid span <.4mm>
									\vshade  0 1 1.9 <,z,,> 1 1 1.9 <z,z,,>  2 0 1.9 <z,z,,> 3 1 1.9 <z,z,,> 7 1 1.9
									  <z,z,,> 8 0 1.9  <z,z,,> 9 1 1.9 <z,,,> 12 1 1.9 /

									  \setsolid
									  \arr{0.3 1.7}{0.7 1.3}
									  \arr{1.3 1.3}{1.7 1.7}
									  \arr{2.3 1.7}{2.7 1.3}
									  \arr{3.3 1.3}{3.7 1.7}
									  \arr{4.3 1.7}{4.7 1.3}
									  \arr{5.3 1.3}{5.7 1.7}
									  \arr{6.3 1.7}{6.7 1.3}
									  \arr{7.3 1.3}{7.7 1.7}
									  \arr{8.3 1.7}{8.7 1.3}
									  \arr{9.3 1.3}{9.7 1.7}
									  \arr{10.3 1.7}{10.7 1.3}
									  \arr{11.3 1.3}{11.7 1.7}
									  \setdashes <1mm>
									  \plot 0 0.3  0 2.3 /
									  \plot 12 0.3  12 2.3 /
									  \multiput{$\cdots$} at 1 2.3  6 2.3  11 2.3 /

									  \endpicture}
									  $$
									  Next, for $n\ge 4$ even and $i =\frac n2$, the pair $([i],[n])$ is the
									  only object in the $\tau_{\Cal S}$-orbit of $([i],[i])$ which survives:
									  $$
									  
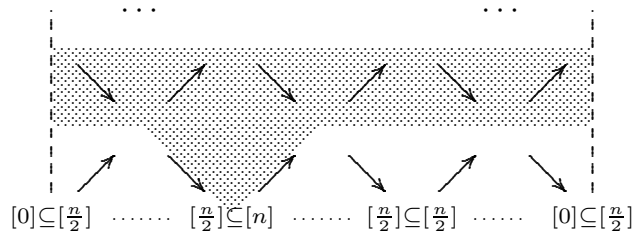
\hbox{\beginpicture
									  \setcoordinatesystem units <1cm,1cm>
									  \put{} at 0 .7

									  \arr{0.3 0.3}{0.7 0.7}
									  \arr{2.3 0.3}{2.7 0.7}
									  \arr{4.3 0.3}{4.7 0.7}
									  \arr{1.3 0.7}{1.7 0.3}
									  \arr{3.3 0.7}{3.7 0.3}
									  \arr{5.3 0.7}{5.7 0.3}

									  \setdots <1mm>
									  \plot .7 0  1.3 0 /
									  \plot 2.7 0  3.3 0 /
									  \plot 4.7 0  5.2 0 /

									  \put{$\ssize[\frac n2]\subseteq[\frac n2]$} at 4 0
									  \put{$\ssize[\frac n2]\subseteq[n]$} at 2 0
									  \multiput{$\ssize[0]\subseteq[\frac n2]$} at 0 0  6 0 /

									  \put{} at 0 2.3

									  \setshadegrid span <.4mm>

									  \vshade  0 1 1.9 <,z,,> 1 1 1.9 <z,z,,>  2 0 1.9 <z,z,,> 3 1 1.9 <z,,,> 6 1 1.9 /
									  \setsolid
									  \arr{0.3 1.7}{0.7 1.3}
									  \arr{1.3 1.3}{1.7 1.7}
									  \arr{2.3 1.7}{2.7 1.3}
									  \arr{3.3 1.3}{3.7 1.7}
									  \arr{4.3 1.7}{4.7 1.3}
									  \arr{5.3 1.3}{5.7 1.7}
									  \setdashes <1mm>
									  \plot 0 0.3  0 2.3 /
									  \plot 6 0.3  6 2.3 /
									  \multiput{$\cdots$} at 1 2.3  5 2.3  /

									  \endpicture}
									  $$
									  Finally, for $i=1$, again two objects in the $\tau_{\Cal S}$-orbit
									  of $([i],[i])$ survive:
									  $$
									  
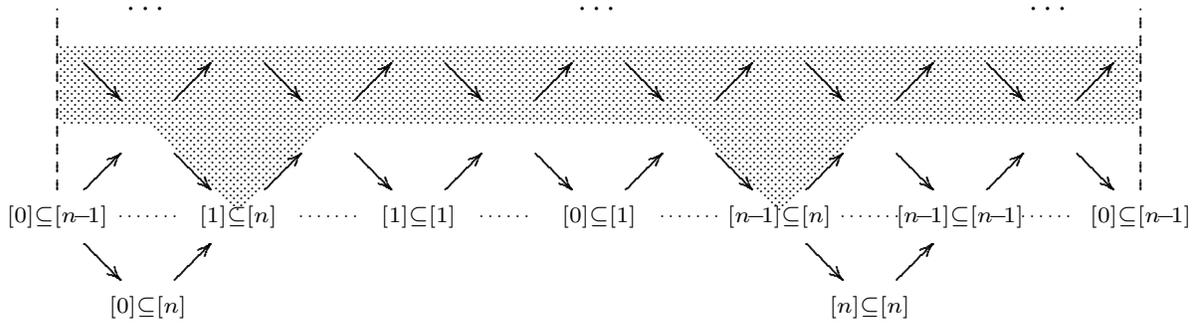
\hbox{\beginpicture
									  \setcoordinatesystem units <1cm,1cm>
									  \put{} at 0 .7
									  \multiput{} at -1 0  13 0 /

									  \put{$\ssize[1]\subseteq[1]$} at 4 0
									  \put{$\ssize[1]\subseteq[n]$} at 2 0
									  \multiput{$\ssize[0]\subseteq[n\!-\!1]$} at 0 0  12 0 /
									  \put{$\ssize[0]\subseteq[n]$} at 1 -1
									  \arr{0.3 0.3}{0.7 0.7}
									  \arr{2.3 0.3}{2.7 0.7}
									  \arr{4.3 0.3}{4.7 0.7}
									  \arr{1.3 0.7}{1.7 0.3}
									  \arr{3.3 0.7}{3.7 0.3}
									  \arr{5.3 0.7}{5.7 0.3}

									  \arr{0.3 -.3}{0.7 -.7}
									  \arr{1.3 -.7}{1.7 -.3}

									  \setdots <1mm>

									  \plot .7 0  1.3 0 /
									  \plot 2.7 0  3.3 0 /
									  \plot 4.7 0  5.2 0 /

									  \setsolid
									  \put{$\ssize[n\!-\!1]\subseteq[n\!-\!1]$} at 10 0
									  \put{$\ssize[n\!-\!1]\subseteq[n]$} at 8 0
									  \put{$\ssize[0]\subseteq[1]$} at 6 0
									  \put{$\ssize[n]\subseteq[n]$} at 9 -1
									  \arr{6.3 0.3}{6.7 0.7}
									  \arr{8.3 0.3}{8.7 0.7}
									  \arr{10.3 0.3}{10.7 0.7}
									  \arr{7.3 0.7}{7.7 0.3}
									  \arr{9.3 0.7}{9.7 0.3}
									  \arr{11.3 0.7}{11.7 0.3}

									  \arr{8.3 -.3}{8.7 -.7}
									  \arr{9.3 -.7}{9.7 -.3}
									  \setdots <1mm>
									  \plot 6.7 0  7.3 0 /
									  \plot 8.7 0  9.3 0 /
									  \plot 10.7 0  11.2 0 /

									  \put{} at 0 2.3

									  \setshadegrid span <.4mm>
									  \vshade  0 1 1.9 <,z,,> 1 1 1.9 <z,z,,>  2 0 1.9 <z,z,,> 3 1 1.9 <z,z,,> 7 1 1.9
									    <z,z,,> 8 0 1.9  <z,z,,> 9 1 1.9 <z,,,> 12 1 1.9 /

									    \setsolid
									    \arr{0.3 1.7}{0.7 1.3}
									    \arr{1.3 1.3}{1.7 1.7}
									    \arr{2.3 1.7}{2.7 1.3}
									    \arr{3.3 1.3}{3.7 1.7}
									    \arr{4.3 1.7}{4.7 1.3}
									    \arr{5.3 1.3}{5.7 1.7}
									    \arr{6.3 1.7}{6.7 1.3}
									    \arr{7.3 1.3}{7.7 1.7}
									    \arr{8.3 1.7}{8.7 1.3}
									    \arr{9.3 1.3}{9.7 1.7}
									    \arr{10.3 1.7}{10.7 1.3}
									    \arr{11.3 1.3}{11.7 1.7}
									    \setdashes <1mm>
									    \plot 0 0.3  0 2.3 /
									    \plot 12 0.3  12 2.3 /
									    \multiput{$\cdots$} at 1 2.3  6 2.3  11 2.3 /

									    \endpicture}
									    $$
										\bigskip
										{\bf Bookkeeping 1.} Going from $\Cal S(n)$ to $\Cal F(n)$, or from $\Cal F(n)$ to
										$\mod\Pi_{n-1}$, the number of indecomposable objects decreases in both step by $n$.
										Here is the bookkeeping table. We denote the number of isomorphism classes of indecomposable objects in the category $\Cal C$ by $\#\ind\Cal C$.
										$$
										\hbox{\beginpicture
										\setcoordinatesystem units <3cm,.4cm>

										\put{$n$} at 0 7.2
										\put{$1$} at 0 6
										\put{$2$} at 0 5
										\put{$3$} at 0 4
										\put{$4$} at 0 3
										\put{$5$} at 0 2
										\put{$6$} at 0 1

										\put{$\#\ind\Cal S(n)$} at 1 7.2
										\put{$2$} at 1 6
										\put{$5$} at 1 5
										\put{$10$} at 1 4
										\put{$20$} at 1 3
										\put{$50$} at 1 2
										\put{$\infty$} at 1 1

										\put{$\#\ind\Cal F(n)$} at 2 7.2
										\put{$1$} at 2 6
										\put{$3$} at 2 5
										\put{$7$} at 2 4
										\put{$16$} at 2 3
										\put{$45$} at 2 2
										\put{$\infty$} at 2 1

										\put{$\#\ind\mod\Pi_{n-1}$} at 3 7.2
										\put{$0$} at 3 6
										\put{$1$} at 3 5
										\put{$4$} at 3 4
										\put{$12$} at 3 3
										\put{$40$} at 3 2
										\put{$\infty$} at 3 1

										\plot -.3 6.7  3.5 6.7 /
										\plot 0.4 8  0.4 0.5 /
										\endpicture}
										$$
												\bigskip 
{\bf Bookkeeping 2.}
Going from $\underline{\Cal S}(n)$ to $\mod \Pi_{n-1}$
the number of indecomposables decreases by $2(n-1)$, going from
$\mod \Pi_{n-1}$ to $\underline{\mod}\ \Pi_{n-1}$ the number decreases
by $n-1$. Here are the actual numbers; for the triangulated categories
$\underline{\Cal S}(n)$ and $\underline{\mod}\ \Pi_{n-1}$ we also
list the tree type of the corresponding Auslander-Reiten quivers.
$$
\hbox{\beginpicture
\setcoordinatesystem units <3cm,.5cm>

\put{$n$} at 0 7.2
\put{$1$} at 0 6
\put{$2$} at 0 5
\put{$3$} at 0 4
\put{$4$} at 0 3
\put{$5$} at 0 2
\put{$6$} at 0 1

\put{$\#\ind\underline{\Cal S}(n)$} at 1 7.2
\put{$0$} at 1 6
\put{$3$} at 1 5
\put{$8$} at 1 4
\put{$18$} at 1 3
\put{$48$} at 1 2
\put{$\infty$} at 1 1
\put{$\#\ind\mod\Pi_{n-1}$} at 2 7.2
\put{$0$} at 2 6
\put{$1$} at 2 5
\put{$4$} at 2 4
\put{$12$} at 2 3
\put{$40$} at 2 2
\put{$\infty$} at 2 1

\put{$\#\ind\underline{\mod}\ \Pi_{n-1}$} at 3 7.2
\put{$0$} at 3 6
\put{$0$} at 3 5
\put{$2$} at 3 4
\put{$9$} at 3 3
\put{$36$} at 3 2
\put{$\infty$} at 3 1

\put{$\Bbb A_2$} at 1.3 5
\put{$\Bbb D_4$} at 1.3 4
\put{$\Bbb E_6$} at 1.3 3
\put{$\Bbb E_8$} at 1.3 2

\put{$\Bbb A_1$} at 3.27 4
\put{$\Bbb A_3$} at 3.27 3
\put{$\Bbb D_6$} at 3.27 2

\plot -.3 6.7  3.5 6.7 /
\plot 0.4 8  0.4 0.5 /
\endpicture}
$$

{\bf 11. Appendix: Objective functors.}
     \medskip
Let $\Cal A, \Cal B$ be additive categories and let $F\:\Cal A \to \Cal B$ be an (additive) functor. 
An object $A$ in $\Cal A$ will be called
a {\it kernel object} for $F$ provided $F(A) = 0.$ 
The functor $F\:\Cal A \to \Cal B$ 
will be said to be {\it objective} provided any morphism $f\: A \to A'$ in $\Cal A$
with $F(f) = 0$ factors through a kernel object for $F$. If $F$ is an objective functor, then
we will say that the {\it kernel of $F$ is generated by} $\Cal K,$ provided 
$\Cal K$ is a class of objects in $\Cal A$ such that $\add \Cal K$ is the 
class of all kernel objects for $F$.

Given an additive category $\Cal A$ and an ideal $\Cal I$ in $\Cal A$,
we denote by $\Cal A/\Cal I$ the corresponding factor category (it has the same
objects, and the homomorphisms in $\Cal A/\Cal I$ are the residue classes of
the homomorphisms in $\Cal A$ modulo $\Cal I$).
If $\Cal K$ is a class of objects of the category $\Cal A$,
we denote by $\langle \Cal K\rangle$ the ideal of $\Cal A$ given by all maps
which factor through a direct sum of objects in $\Cal K.$ Instead of writing
$\Cal A/\langle \Cal K\rangle$, we just will write $\Cal A/\Cal K.$

If $F\:\Cal A \to \Cal B$ is a full, dense, objective functor and the kernel of $F$ 
is generated by $\Cal K$, then $F$ induces an equivalence between the category $\Cal A/\Cal K$ and $\Cal B$.
We see that given a full, dense, objective functor $F\:\Cal A \to \Cal B$, the category $\Cal B$ 
is uniquely determined by $\Cal A$ and a class of indecomposable objects in $\Cal A$ (namely the class 
of indecomposable kernel objects for $F$); if $F$ is objective, but not necessarily full or dense, then
$F$ induces an equivalence between the category $\Cal A/\Cal K$ and the image category of $F$.
     \medskip
{\it The composition of objective functors is not necessarily objective.} Here is an 
example: Let $\Cal B$ be the linearization of the chain of cardinality 3, thus $\Cal B$ has three
objects $b_1,b_2,b_3$ with $\Hom(b_i,b_j) = k$ provided $i \le j$ and zero otherwise, such that
the composition $\Hom(b_2,b_3)\otimes \Hom(b_1,b_2) \to \Hom(b_1,b_3)$ is the multiplication map. 
Let $\Cal A$ be the full subcategory of $\Cal B$ with objects $b_1,b_3$, thus $\Cal A$ is the
linearization of a chain of cardinality 2. Let $\Cal K = \{b_2\}$ and $\Cal C = \Cal B/\Cal K$.
The inclusion functor $F\:\Cal A \to \Cal B$ and the projection functor $G\:\Cal B \to \Cal C$ both
are (full and) objective, however the composition $GF\:\Cal A \to \Cal B$ is not objective (none of the objects
$b_1,b_3$ belongs to the kernel of $GF$, we have $\Hom_{\Cal A}(b_1,b_3) = k$ and 
any non-zero map $b_1 \to b_3$ is mapped to zero under $GF$). Note that the functor $F$ is not dense.
    \medskip
{\it  Let $F\:\Cal A \to \Cal B$ and $G\:\Cal B \to \Cal C$ be objective functors. If $F$ is, in addition, full and 
dense, then $GF$ is objective}  (thus, the composition of full, dense, objective functors is full, dense, objective).
Proof. Since $F, G$ both are full, also $GF$ is full. Let $a\:A_1\to A_2$ be a morphism with $GF(a) = 0.$
Since $G$ is objective, the morphism $F(a)$ factors through a kernel object $B$ for $G$, say $F(a) = b_2b_1$
where $b_1\:F(A_1) \to B$ and $b_2\:B \to F(A_2)$. By assumption, $F$ is dense, thus there is an
isomorphism $b\:B \to F(A)$ for some object $A$ in $\Cal A$. Since $F$ is full, there is a map $a_1\:A_1 \to A$
such that $F(a_1) = bb_1$ and a map $a_2\:A \to A_2$ 
such that $F(a_2) = b_2b^{-1}$. We have $F(a) = b_2b^{-1}bb_1 = F(a_2)F(a_1) = F(a_2a_1),$  thus 
$F(a-a_2a_1) = 0.$ Since $F$ is objective, there is a kernel object $A'$ for $F$ such that $a-a_2a_1$
factors through $A'$, say $a-a_2a_1 = a_4a_3$, with $a_3\:A_1\to A', a_4\:A' \to A_2.$ It follows that
$a = a_2a_1+a_4a_3 = \bmatrix a_2 &a_4\endbmatrix \bmatrix a_1 \cr a_3\endbmatrix,$ thus this map factors
through $A\oplus A'$. But $GF(A\oplus A') = GF(A)\oplus GF(A').$ Now, $F(A)$ is isomorphic to $B$, thus
$GF(A)$ is isomorphic to $G(B) = 0$. Also, $F(A') = 0,$ thus $GF(A') = 0$. This shows that $A\oplus A'$ is
a kernel object for $GF.$
	\medskip
We recall that an additive category $\Cal A$ is said to be
a Krull-Remak-Schmidt category, provided every object in $\Cal A$ is a (finite) direct sum of
objects with local endomorphism rings. Assume now that  $F\: \Cal A \to \Cal B$
is an objective functor between Krull-Remak-Schmidt categories $\Cal A$ and $\Cal B$.
Then we are interested in the number $i_0(F)$ of isomorphism classes of indecomposable
objects in $\Cal F$ which are kernel objects for $F$, 
as well as in the number $i_1(F)$ of isomorphism classes
of indecomposable objects $B$ in $\Cal B$ such that $B$ is not isomorphic to an object of the form $F(A)$
where $A$ is an object in $\Cal A$. If at least one of the numbers $i_0(F), i_1(F)$ is finite, we call
$i(F) = i_0(F)-i_1(F)$ the {\it index} of $F$. 

The objective functors $F$ considered in the paper are also dense, 
in this case $i(F) = i_0(F)$ is the number of
isomorphism classes of indecomposable kernel objects in $\Cal A$.
  \bigskip\bigskip 

  {\bf 12. References.}
       \medskip
       \item{[A]} A.~Auslander: Coherent functors. In: Proceedings of the Conference on
         Categorical Algebra. La Jolla 1965, Springer-Veriag, New York (1965), 189--231.
	 \item{[AR1]} M.~Auslander, I.~Reiten: Stable equivalence of dualizing $R$-varieties.
	    Adv. Math. 12 (1974), 306--366.
	    \item{[AR2]} M.~Auslander, I.~Reiten: On the representation type of triangular
	      matrix rings.  J. London Math. Soc.(2), 12 (1976), 371--382.
	      \item{[BMR]} A.~B.~Buan, R.~J.~Marsh, I.~Reiten: Cluster-tilted algebras. Trans. Amer. Math.
	        Soc. 359 (2007), 323--332.
		\item{[DR]} V.~Dlab, C.~M.~Ringel: The module theoretical approach to
		 quasi-hereditary algebras.
		 In: Representations of Algebras and Related Topics (ed. H.~Tachikawa and S.~Brenner). London Math. Soc. Lecture Note Series 168.
		 Cambridge University Press (1992), 200--224.
		 \item{[K]} B.~Keller: On triangulated orbit categories. Doc. Math. 10 (2005), 551--581.
		 \item{[LZ]} Z.-W.~Li, P.~Zhang: A construction of Gorenstein-projective
		  modules. J. Algebra 323 (2010), 1802--1812.
		  \item{[R]}  C.~M.~Ringel:  Iyama's finiteness theorem via strongly quasi-hereditary algebras. Journal Pure Applied Algebra 214 (2010), 1687--1692.
		  \item{[RS1]} C.~M.~Ringel, M.~Schmidmeier: The Auslander-Reiten translation
		   in submodule categories. Trans. Amer. Math. Soc. 360 (2008), 691--716.
		   \item{[RS2]} C.~M.~Ringel, M.~Schmidmeier: Invariant subspaces of nilpotent
		    linear operators. I. J. reine angew. Math 614 (2008), 1--52.
\item{[Z]} P.~Zhang,   Gorenstein-projective modules and symmetric recollements
   J. Algebra 388 (2013), 65--80.
			\bigskip\bigskip
			{\rmk
			C. M. Ringel\par
			Department of Mathematics, \par
			Shanghai Jiao Tong University \par
			Shanghai 200240, P. R. China, \par
			and \par
			King Abdulaziz University, P O Box 80200\par
			Jeddah, Saudi Arabia\par
			e-mail: {\ttk ringel\@math.uni-bielefeld.de} \par
			\medskip
			P. Zhang\par
			Department of Mathematics, \par
			Shanghai Jiao Tong University \par
			Shanghai 200240, P. R. China, \par
			e-mail: {\ttk pzhang\@sjtu.edu.cn} \par
			}
				\bigskip
				\bye